\makeatletter

\def\eqalign#1{\,\vcenter{\openup\jot\m@th
  \ialign{\strut\hfil$\displaystyle{##}$&$\displaystyle{{}##}$\hfil
      \crcr#1\crcr}}\,}
\def\eqalignno#1{\displ@y \tabskip\@centering
  \halign to\displaywidth{\hfil$\displaystyle{##}$\tabskip\z@skip
    &$\displaystyle{{}##}$\hfil\tabskip\@centering
    &\llap{$##$}\tabskip\z@skip\crcr
    #1\crcr}}
\def\leqalignno#1{\displ@y \tabskip\@centering
  \halign to\displaywidth{\hfil$\displaystyle{##}$\tabskip\z@skip
    &$\displaystyle{{}##}$\hfil\tabskip\@centering
    &\kern-\displaywidth\rlap{$##$}\tabskip\displaywidth\crcr
    #1\crcr}}

\makeatother

\newdimen\pixel \pixel=.00333333 true in

\documentclass[11pt,a4paper,reqno]{article}
\usepackage{amssymb,amsmath,amscd,epsfig,enumerate}

\def\footnoterule{\kern-.5\pixel
        \hrule height \pixel width \columnwidth
        \kern-.5\pixel}

\addtocounter{footnote}{1}

\makeatletter
\def\bigpar{\bigbreak\@afterindentfalse\@afterheading\ignorespaces}
\def\medpar{\medbreak\@afterindentfalse\@afterheading\ignorespaces}
\def\smallpar{\smallbreak\@afterindentfalse\@afterheading\ignorespaces}

\newlength{\saveindent}
\newenvironment{proof}%
      {\bigpar{\bf Proof:}\ %  previously \sentsp rather than \ 
             \setlength{\saveindent}{\parindent} %\setlength{\parindent}{0pt}%
                       \ignorespaces}%
      {\stopproof\ignorespaces\bigbreak \setlength{\parindent}{\saveindent}}

\newenvironment{proofnobox}%
      {\bigpar{\bf Proof:}%
             \setlength{\saveindent}{\parindent} %\setlength{\parindent}{0pt}%
                       \,\ignorespaces}%
      {\ignorespaces\bigbreak \setlength{\parindent}{\saveindent}}

      {\bigpar{\bf Proof:}\ %
             \setlength{\saveindent}{\parindent} %\setlength{\parindent}{0pt}%
                       \ignorespaces}%
      {\ignorespaces\bigbreak \setlength{\parindent}{\saveindent}}

\newenvironment{proofof}[1]%
      {\bigpar{\bf Proof of #1:}\ %
             \setlength{\saveindent}{\parindent} %\setlength{\parindent}{0pt}%
                       \ignorespaces}%
      {\stopproof\ignorespaces\bigbreak \setlength{\parindent}{\saveindent}}

\newenvironment{proofofnobox}[1]%
      {\bigpar{\bf Proof of #1:}\ %
             \setlength{\saveindent}{\parindent} %\setlength{\parindent}{0pt}%
                       \ignorespaces}%
      {\ignorespaces\bigbreak \setlength{\parindent}{\saveindent}}

      {\medpar{\bf Remark:}\ %  **** \bigpar, and \bigbreak below ****
                       \ignorespaces\small}%
      {\stopproof\ignorespaces\medbreak \setlength{\parindent}{\saveindent}} 

\newenvironment{remark*}%
      {\medpar{\bf Remark:}\ %  **** \bigpar, and \bigbreak below ****
                       \ignorespaces\small}%
      {\ignorespaces\medbreak \setlength{\parindent}{\saveindent}} 

      {\medpar{\bf Remarks:}\ % **** ditto ****
                       \ignorespaces}%
      {\stopproof\ignorespaces\medbreak \setlength{\parindent}{\saveindent}}

\newenvironment{remarks*}%
      {\medpar{\bf Remarks:}\ % **** ditto ****
                       \ignorespaces}%
      {\ignorespaces\medbreak \setlength{\parindent}{\saveindent}}

      {\medpar{\bf Remark:}\ %
                       \ignorespaces}%
      {\stopproof\ignorespaces\medbreak \setlength{\parindent}{\saveindent}}

      {\medpar{\bf Remarks:}\ %
                       \ignorespaces}%
      {\stopproof\ignorespaces\medbreak \setlength{\parindent}{\saveindent}}
\makeatother

\newtheorem{theorem}{Theorem}[section]    % Remove percent sign
\newtheorem{lemma}[theorem]{Lemma}
\newtheorem{claim}[theorem]{Claim}
\newtheorem{proposition}[theorem]{Proposition}
\newtheorem{corollary}[theorem]{Corollary}
\newtheorem{definition}[theorem]{Definition}

\newtheorem{example}{Example}
\def\begex{\begin{example}\parindent=0pt \rm}
\def\endex{\end{example}}
\def\square{\vbox{\hrule height.2pt\hbox{\vrule width.2pt height5pt \kern5pt
                                   \vrule width.2pt} \hrule height.2pt}}
\def\stopproof{\qquad\square}

\newskip\storeadskip
\newskip\storebdskip
\newcounter{sparectr}

                          {\end{list}}          

                          {\end{list}}          

                          {\end{list}}          

\newcounter{sparectrtwo}

                          {\end{list}}          

\newcounter{sparectrthree}
                          {\end{list}}          
                          {\end{list}}

\newcounter{refcount}

\def\half{{\textstyle{1\over2}}}

\input amssym.def
\def\Rset{{\Bbb R}}

\def\Zset{{\Bbb Z}}

\font\teneusm=eusm10
\font\seveneusm=eusm7
\font\fiveeusm=eusm5
\newfam\eusmfam
\textfont\eusmfam=\teneusm
\scriptfont\eusmfam=\seveneusm
\scriptscriptfont\eusmfam=\fiveeusm

\def\Tree{{\Bbb T}}
\def\muAeta{{\mu_A^\eta}}
\def\VarAeta{{\Var_A^\eta}}
\def\EntAeta{{\Ent_A^\eta}}

\def\const{{\rm const}}
\def\height{{\rm height}}
\def\Expect{{\rm E}}
\def\Plus{$(+)$}
\def\Minus{$(-)$}

\newcommand{\set}[1]{\left\{#1\right\}}

\newcommand{\Var}[0]{{\mathrm{Var}}}
\newcommand{\PV}[0]{{\mathrm{Pvar}}}

\newcommand{\E}[0]{\mu}
\newcommand{\Ent}[0]{{\mathrm{Ent}}}

\newcommand{\Z}[0]{\mathbb{Z}}

\newcommand{\cgap}[0]{c_{\mathrm{gap}}}
\newcommand{\csob}[0]{c_{\mathrm{sob}}}
\newcommand{\df}[0]{{\mathcal{D}}}
\newcommand{\dfe}[0]{{\mathcal{E}}}
\newcommand{\vcond}[0]{{\hbox{\rm VM}}}
\newcommand{\econd}[0]{{\hbox{\rm EM}}}
\newcommand{\Cov}[0]{{\mathrm{Cov}}}
\newcommand{\partialb}[0]{\widetilde{\partial}}

 \let\b=\beta   \let\d=\delta  
  \let\h=\eta      
            
  \let\s=\sigma \let\t=\tau   \let\th=\vartheta

\let\O=\Omega

\newcommand{\cF}{\ensuremath{\mathcal F}}

\newcommand{\cL}{\ensuremath{\mathcal L}}

\let\neper=e
\let\ii=i

\def\ie{\hbox{\it i.e.\ }}

\def\nep#1{ \neper^{#1}}

\def\ov#1{{1\over#1}}

\def\tc{\thsp | \thsp}
\let\<=\langle
\let\>=\rangle

\def\Ent{ \mathop{\rm Ent}\nolimits }
\def\Cov{ \mathop{\rm Cov}\nolimits }

\def\ninf#1{ \| #1 \|_\infty }

\outer\def\nproclaim#1 [#2]#3. #4\par{\medbreak \noindent
   \talato(#2){\bf #1 \Thm[#2]#3.\enspace }%
   {\sl #4\par }\ifdim \lastskip <\medskipamount
   \removelastskip \penalty 55\medskip \fi}
\def\thmm[#1]{#1}
\def\teo[#1]{#1}
\def\sttilde#1{%
\dimen2=\fontdimen5\textfont0
\setbox0=\hbox{$\mathchar"7E$}
\setbox1=\hbox{$\scriptstyle #1$}
\dimen0=\wd0
\dimen1=\wd1
\advance\dimen1 by -\dimen0
\divide\dimen1 by 2
\vbox{\offinterlineskip%
   \moveright\dimen1 \box0 \kern - \dimen2\box1}
}

\def\smallno{\smallskip\noindent}
\def\medno{\medskip\noindent}

\def\\{\hfill\break}

\def\thsp{\thinspace}

\def\tthsp{\kern .083333 em}

\def\?{\mskip -10mu}

\renewcommand{\include}{\input}
\begin{document}
% \include{title_CMP}
% \include{intro_CMP}
% \include{acknowledgments_CMP}
% \include{prelims_CMP}
% \include{spat_CMP}
% \include{gap_CMP}
% \include{logsob_CMP}
% \include{extensions_CMP}
% \include{supplement_CMP}
% \include{refs_CMP}
%%%%%%%%%%%%%%%%%% TITLE PAGE %%%%%%%%%%%%%%%%%%
\title{Glauber dynamics on trees:\\ 
Boundary conditions and mixing time\thanks{An extended abstract
of this paper appeared under the title ``The Ising model on trees:
Boundary conditions and mixing time'' in {\it Proceedings of the
44th Annual IEEE Symposium on Foundations of Computer Science},
October 2003.}
%\\(Preliminary Version: Please Do Not Distribute!)
}

\author{Fabio Martinelli\thanks{Department of Mathematics, University
of Roma Tre, Largo San Murialdo~1, 00146~Roma, Italy.\ \ Email:
{\tt martin@mat.uniroma3.it}.  This work was done while this author
was visiting the Departments of EECS and Statistics,
University of California, Berkeley, supported in
part by a Miller Visiting Professorship.}\and
Alistair Sinclair\thanks{Computer Science Division, University
of California, Berkeley, CA~94720-1776, U.S.A.\ \ Email:
{\tt sinclair@cs.berkeley.edu}.  Supported in part by
NSF Grant CCR-0121555 
and DARPA cooperative agreement F30602-00-2-0601.}\and
Dror Weitz\thanks{Computer Science Division, University
of California, Berkeley, CA~94720-1776, U.S.A.\ \ Email:
{\tt dror@cs.berkeley.edu}.  
Supported in part by NSF Grant CCR-0121555.}}
\maketitle
\thispagestyle{empty}
%%%% abstract %%%%%%
\begin{center}
\large\bf Abstract
\end{center}
%\begin{abstract}
\noindent 
We give the first comprehensive analysis of the effect of boundary
conditions on the mixing time of the Glauber dynamics in the so-called
\emph{Bethe approximation}.  Specifically, we show that spectral gap
and the log-Sobolev
constant of the Glauber dynamics for the Ising model on an $n$-vertex 
regular tree with \Plus-boundary are bounded below by a constant independent
of~$n$ at all temperatures and all external fields.  This implies that
the mixing time is $O(\log n)$
(in contrast to the free boundary case, where it is not
bounded by any fixed polynomial at low temperatures).
In addition, our methods yield
simpler proofs and stronger results for the spectral gap and 
log-Sobolev constant in the regime where there are multiple phases but the
mixing time is insensitive to the boundary condition. 
Our techniques also apply to a much wider class of models, including 
those with hard-core
constraints like the antiferromagnetic Potts model at zero temperature
(proper colorings) and the hard--core lattice gas (independent sets).
%\end{abstract}
\newpage
\setcounter{footnote}{1}
\setcounter{page}{1}

%%%%%%%%%%%%%%%%%% INTRODUCTION %%%%%%%%%%%%%%%%%%
\section{Introduction}
\label{sec:intro}
% \vskip-0.3in\hbox{}
% \subsection{Background}
% \vskip-0.05in
In this paper we will analyze the influence of boundary conditions on
the Glauber dynamics for discrete spin models on a regular rooted tree. 
Although in what follows we will focus for simplicity on the well known 
Ising  model, 
our techniques also apply to other models, not necessarily
ferromagnetic and with hard-core constraints.

In the Ising model on a finite graph $G=(V,E)$, a configuration
$\sigma=(\sigma_x)$
consists of an assignment of $\pm 1$-values, or ``spins'', to
each vertex (or ``site'') of~$V$.  The probability of finding the system in
configuration $\sigma\in\{\pm 1\}^V\equiv\Omega_G$ is given by the 
{\it Gibbs distribution\/}
\begin{equation}\label{eq:gibbs}
  \mu_G(\sigma) \propto \exp\Bigl(\beta
\sum\nolimits_{xy\in E}\sigma_x\sigma_y+\b h \sum\nolimits_{x\in V}\sigma_x\Bigr),
\end{equation}
where $\beta\ge 0$ is the inverse temperature and $h$ the external field.
Boundary conditions can also be taken into account by fixing the
spin values at some specified ``boundary'' vertices of~$G$; the term
{\it free boundary\/} is used to indicate that no boundary condition
is specified.

In the classical Ising model, $G=G_n$ is a cube of side~$n^{1/d}$
in the $d$-dimensional Cartesian lattice~$\Zset^d$, and in this case the phase
diagram in the thermodynamic limit $G_n\uparrow \Z^d$ 
is quite well understood (see, e.g., \cite{Georgii,Simon} for more
background).

While the classical theory focused on static properties of the Gibbs
measure, in the last decade the emphasis has shifted
towards dynamical questions with a computational flavor. The key 
object here is the {\it Glauber dynamics}, a (discrete-- or continuous--time) Markov chain on the
set of spin configurations~$\Omega_G$ in which each spin $\s_x$ flips
its value with a rate that depends on the current configuration of the
neighboring spins of $x$, and which satisfy the detailed balance
condition w.r.t to the Gibbs measure $\mu_G$ (see Section~\ref{sec:prelims} 
for more details).

The Glauber
dynamics is much studied for two reasons: firstly, it is the
basis of Markov chain Monte Carlo algorithms, widely used in
computational physics for sampling from the Gibbs distribution;
and secondly, it is a plausible model for the actual evolution
of the underlying physical system towards equilibrium. In both
contexts, one of the central questions is to determine the {\it mixing time},
i.e., the time until the dynamics is close to its stationary
distribution. 

As is well known (see e.g. \cite{Saloff}), the approach to stationarity
of a reversible Markov chain with Markov generator $\cL$ and reversible
measure $\pi$ can be succesfully studied by analyzing two key
quantities: the \emph{spectral gap} and the \emph{logarithmic Sobolev
  constant} of the pair $(\cL,\pi)$\footnote{Unfortunately the
  definition of the logarithmic Sobolev constant is not constant in the
  literature. The ambiguity arises because there are two definitions,
  one the inverse of the other. The definition used in this paper
  is the one that puts the logarithmic Sobolev constant and the
  spectral gap on the same footing.}.  The first of these measures the
rate of the exponential decay as $t\to \infty$ of the variance
$\Var_\pi(\nep{t\cL}f)$ computed with respect to the invariant measure
$\pi$, while the second measures instead the rate of decay of the
relative entropy of $\nep{t\cL}f$ w.r.t $\pi$ (see, e.g., \cite{bible}).
Advances in statistical physics over the past decade have led to
remarkable connections between these two quantities and the occurence of
a phase transition (see, e.g., \cite{SZ,MOS,MO1,Cesi,Martinelli,Yau}).
As an example, on finite $n$-vertex squares with free boundary in the
2-dimensional lattice~$\Zset^2$, when $h=0$ and $\beta$ is smaller than
the critical value $\beta_c$, the spectral gap and the logarithmic
Sobolev constant are $\O(1)$ (\ie bounded away from zero uniformly in
$n$), while for $\beta>\beta_c$ they are both exponentially small in
$\sqrt{n}$.

One of the most interesting and difficult questions left open by the
above and related results is the {\it influence of boundary
conditions\/} on the spectral gap and the log-Sobolev constant 
when $h=0$ and $\b>\b_c$.  It
has been conjectured that, in the presence of an all-\Plus\ boundary, the
relaxation process is driven by the mean--curvature motion of interfaces separating
droplets of the \Minus-phase inside the \Plus-phase, and therefore the
mixing time should be polynomial in~$n$ (most likely $n^{2/d}\log n$)
~\cite{BM,FH}.  In particular it has been argued that the
spectral gap for the pure phases in high enough dimension should be
$\O(1)$. Proving results of this kind has proved very elusive, and
the only (presumably sharp) available bounds are
\emph{upper bounds} on the spectral gap and the logarithmic Sobolev
constant \cite{BM}.

In this paper we prove a strong version of the above conjecture in what
is known in statistical physics as the {\it Bethe approximation}, namely
when the lattice $\Zset^d$ is replaced by a regular tree.  Among other
results, we show that the spectral gap of the Glauber dynamics for the
Ising model on a tree with a \Plus-boundary condition on its leaves is
$\O(1)$ at all temperatures and all values of the external field, and 
further that the same holds for the logarithmic Sobolev constant.
Notice that, with a free boundary, $\b$ large and $h=0$, 
both quantities tend to zero as $1/n^a$ and the exponent
$a$ grows arbitrarily large as $\beta\to\infty$~\cite{BKMP}.

Ours is apparently the first result that quantifies the effect of
boundary conditions on Glauber dynamics in an interesting scenario.  We
stress that, while the tree is simpler in many respects than~$\Zset^d$
due to the lack of cycles, in other respects it is more complex due to
the large boundary: e.g.,
it exhibits a ``double phase transition,'' and the critical field at low
temperature is non--zero (see below). In the next
subsection, we briefly describe the Ising model on trees before stating
our results in more detail.

%\vskip-0.3in\hbox{}
\subsection{The Ising model on trees}
\label{sec:Ising_on_trees}
\begin{figure}[h]
\centerline{\psfig{file=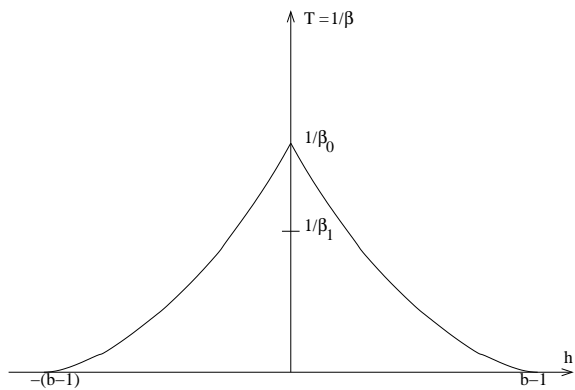,height=2in}}
\caption{The critical field $h_c(\beta)$.  The Gibbs
measure is unique above the curve.}
\label{fig:phase_diagram}
\end{figure}

%\vskip-0.05in 
Fix $b\ge 2$ and let $\Tree^b$ denote the infinite $b$-ary
tree.  The Ising model on~$\Tree^b$ is known \cite{Georgii,L} to have a
phase diagram in the $(h,\b)$ plane quite different from that on the
cubic lattice $\Z^d$ (see Fig.~\ref{fig:phase_diagram}), and has recently
received a lot of attention as the canonical example of a statistical
physics model on a ``non-amenable'' graph (i.e., one whose boundary is
of comparable size to its volume) --- see, e.g.,
\cite{BRZ,Ioffe,EKPS,ST98,JS99,BKMP,BRSSZ}.

Let us first discuss the behavior on the line $h=0$.  There is
a first critical value \hbox{$\b_0= \ov2
\log\bigl(\frac{b+1}{b-1}\bigr)$}, marking the dividing line between uniqueness
and non-uniqueness of the Gibbs measure. Then, in sharp contrast to the
model on~$\Zset^d$, there is a second critical point 
\hbox{$\b_1= \ov2\log\bigl(\frac{\sqrt{b}+1}{\sqrt{b}-1}\bigr)$}
which is often
referred to as the ``spin-glass critical point''~\cite{CCST}. This
second critical point is such that, in the ``intermediate temperature''
region $\beta_0<\beta\le\beta_1$, the \Plus- and \Minus-boundary
conditions exert arbitrarily long-range influence on the spin at root of
the tree and hence give rise to different Gibbs measures, but
``typical'' boundary conditions (i.e., chosen from the infinite
volume Gibbs measure with free boundary) do not. Another way to phrase
this peculiar behavior is that the Gibbs measure constructed via a free
boundary is {\it extremal\/} for all $\b\le \b_1$ (see
~\cite{BRZ,Ioffe,Ioffe2,BKMP} and also \cite{EKPS,Elchsurvey,MP01} for
an analysis in the context of ``bit reconstruction problems'' for noisy
data transmission).

Let us now examine what happens when an external field $h$ is added to
the system. It turns out that for all $\beta>\beta_0$, there is a
critical value~$h=h_c(\beta)>0$ of the field such that the Gibbs measure
is not unique when $|h|\le h_c$, and is unique when $|h|>h_c$.  (When
$\beta\le\beta_0$ the Gibbs measure is unique for all~$h$, and $h_c$ is
defined to be zero.) In the presence of
a \Plus-boundary, the Ising model on the tree with external field
$h=-h_c$ is rather analogous to the classical case of~$\Zset^d$ with
zero field. Both models share the following two properties: firstly, the
Gibbs measure is sensitive to the choice of boundary condition, 
and secondly, adding an arbitrarily small negative field causes the 
Gibbs measure to become insensitive
to the boundary condition (i.e., unique in the thermodynamic limit).

Finally we remark that the concentration properties of the Gibbs measure
for $\b>\b_0$, $h\ge -h_c$ and \Plus-boundary are very different 
from those on~$\Z^d$. 
In the latter case, along the line of first order phase
transition, the (negative) large deviations for the bulk magnetization 
are related
to the appearance of a Wulff droplet of the opposite phase and are
depressed by a negative exponential in the \emph{surface} of the droplet
(see, e.g.,~\cite{DKS}). Here instead, for any value of $(\b,h)$ they are
always depressed by a negative exponential in the \emph{volume} of
the excess negative spins (the phenomenon of 
``rigidity of the critical phases''~\cite{BRSSZ}).
  
The Glauber dynamics for the Ising model on trees has also been
studied.  In a recent paper~\cite{BKMP}, it is shown that the associated
spectral gap (see (\ref{eq:cgapcsob}) for a precise definition) with
zero external field and free boundary on a complete $b$-ary
tree $T$ with $n$~vertices is $\O(1)$ at high and intermediate
temperatures (i.e., when $\beta<\beta_1$)\footnote{Actually the
arguments in \cite{BKMP} prove that the gap is $\O(1)$ for any $\b<\b_1$,
arbitrary boundary condition and any external field. Their
argument, together with some monotonicity properties specific to the 
Ising model~\cite{PW}, implies a mixing time of $O(\log n)$. 
Thus, although for $\b_0<\b<\b_1$ there exist several Gibbs
measures, the mixing time of the Glauber dynamics is
\emph{insensitive} to the boundary condition.}. Moreover, at the critical
point $\b=\b_1$ the same spectral gap is bounded above by $c/\log n$, and
as soon as $\beta>\beta_1$ it becomes smaller than $c/n^{a(\b)}$, with
$a(\b)\uparrow \infty$ as $\b\to \infty$. Thus the critical point
$\beta=\beta_1$ is reflected in the dynamics by an abrubt jump in the
behavior of the spectral gap as a function of the size of the tree $T$.
Finally, also in \cite{BKMP}, it is proved that the spectral gap for
arbitrary fixed $\b,h$ and boundary condition can never shrink to zero
faster than an inverse polynomial in $n$. Again such a result should be
compared to the lattice case where it is known that the spectral gap for a
cube with $n$ sites can be exponentially small in the surface
$n^{(d-1)/d}$.

%\vskip-0.3in\hbox{}
\subsection{Main results and techniques}
%\vskip-0.05in 
Our first main result is a detailed analysis of the spectral
gap of the Glauber dynamics in different regions of the phase diagram. 
The main
novelty here is that we are able for the first time to prove a sharp
result in the region where the spectral gap is highly sensitive to the
boundary condition. 
\begin{theorem}
\label{A}
In both of the following situations, the spectral gap of the Glauber
dynamics on a
complete $b$-ary tree $T$ with $n$~vertices is $\O(1)$:
\begin{itemize}
\item[{\rm (i)}] the boundary condition is arbitrary, 
and either $\beta<\beta_1$ (with $h$
  arbitrary), or $|h|>h_c(\beta)$ (with $\beta$ arbitrary);
\item[{\rm (ii)}] the boundary condition is~\Plus\  and
  $\beta,h$ are arbitrary.
\end{itemize}
\end{theorem}
\begin{remark*}
On $\Z^d$ not much is known about the spectral gap when $\b>\b_c$,
$h=0$ and the boundary condition is~\Plus, the notable exception being
that of $\Z^2$ where it has been recently proved~\cite{BM} that the
spectral gap in a square with $n$ sites shrinks to zero at least as
$1/\sqrt{n}$ (neglecting logarithmic corrections).  The best known lower bounds are significantly 
weaker~\cite{Martinelli}.
In high enough dimensions ($d\ge 3$) it has been conjectured
(see \cite{FH} and \cite{BM}) that the spectral gap should stay bounded
away from zero uniformly in~$n$. The above theorem can be looked upon as 
evidence in favor of this conjecture.  
\end{remark*}
In our second main result we extend our analysis to the more
delicate and difficult logarithmic Sobolev constant (see
(\ref{eq:cgapcsob}) for a precise definition).
\begin{theorem}
\label{B}
In the same situations as in Theorem~\ref{A}, the logarithmic Sobolev constant
of the Glauber dynamics on a complete $b$-ary tree $T$ with $n$~vertices 
is $\O(1)$.
%\begin{itemize}
%\item[{\rm (i)}] the boundary condition is arbitrary, 
%and either $\beta<\beta_1$ (with $h$
%  arbitrary), or $|h|>h_c(\beta)$ (with $\beta$ arbitrary);
%\item[{\rm (ii)}] the boundary condition is~\Plus, $\beta$ is arbitrary
%  and $h\neq -h_c(\b)$.
%\end{itemize}
\end{theorem}
As a corollary we obtain that, in the situations of Theorems~\ref{A} 
and~\ref{B}, the Glauber dynamics mixes (in a very strong sense) in 
time $O(\log n)$.
\begin{remarks*}
  \begin{enumerate}[(i)]
  \item In $\Z^d$ with \Plus-boundary condition, $\b$ large and zero
    external field the logarithmic Sobolev constant in a cube with $n$
    sites is always smaller than $n^{-2/d}$, neglecting logarithmic
    corrections \cite{BM}, in agreement with heuristic predictions based
    on mean--curvature motion of phases interfaces.
\item
  We also prove (see Theorem~\ref{th:crudelogsob}) an additional 
  result which shows that, for an arbitrary nearest-neighbor spin system
  on a tree, as soon as the spectral gap is
  $\O(1)$ then the logarithmic Sobolev constant cannot shrink faster than
  $(c\log n)^{-1}$.  This means that, even when
  a constant lower bound is known for the gap but not for log-Sobolev,
  one can deduce a mixing time of $O((\log n)^2)$.  
  While we do not require this fact to derive the results of this paper, 
  we believe it may be of interest for other models on trees.
  \end{enumerate}
\end{remarks*}

In order to better appreciate Theorem~\ref{B}, one should keep in mind
that for general finite range, translation invariant, compact spin models
on $\Z^d$, if there exists 
an infinite volume Gibbs measure $\mu$ with a positive logarithmic
Sobolev constant, then the system is necessarily in the uniqueness
region and $\mu$ has 
exponentially decaying correlations \cite{SZ2}\footnote{A close look
at the proof in \cite{SZ2} reveals that the same is true for any
infinite, 
locally finite, bounded degree graph such that the volume of any ball of
radius $\ell$ grows sub--exponentially in $\ell$.}. We also recall 
(see, e.g.,~\cite{Ledoux}) that when the log-Sobolev constant is
bounded away from zero one can derive very strong
(Gaussian--like) concentration properties of the corresponding Gibbs
measure, such as those proved in~\cite{BRSSZ}.

We now proceed to sketch some of our techniques and point out
the main technical innovations.  

Our analysis of both the log-Sobolev constant and the spectral gap rests
on certain {\it spatial mixing conditions\/} that can be stated as
follows. Let $f$ be a function %(non--negative in the case of the entropy)
of the spin configuration that \emph{does not depend} on the spins
in the first $\ell$ levels of the tree starting from the root~$r$, 
and let $\mu(f\tc \s_r)$ be the projection of $f$ onto the spin
$\s_r$ at the root. If the variance (respectively, the entropy) under
the Gibbs measure $\mu$ of $\mu(f\tc \s_r)$ decays fast enough with the
depth $\ell$, then we show by a unified argument how to deduce a bound
of~$\O(1)$ on the spectral gap (respectively, the log-Sobolev constant).
Crucially, in contrast to previous approaches we do not require the
above decay to hold in arbitrary environments, but only for the Gibbs
measure $\mu$ under consideration.   This opens up 
the possibility that the condition holds for some boundary
conditions and not for others (with the same values of temperature
and external field).
We also prove the converse, thus
showing the that our mixing conditions are in fact equivalent to the
required bounds on the spectral gap and log--Sobolev constants.

This analysis has several advantages over previous ones~\cite{BKMP,PW}:
it is more direct, applies also when there is an external field,
and applies to general nearest-neighbor spin systems on trees.

The second main ingredient of the paper is establishing the
above spatial mixing conditions in the scenarios of interest described
in the above two theorems.  This is done via a rather simple and novel
coupling technique for the case of the variance. Such a technique
provides, along the way, a new and really elementary proof of the
extremality of the Gibbs measure with free boundary below~$\b_1$.

Surprisingly, we are also able to exploit the same coupling technique
(via strong concentration properties of the Gibbs measure)
to establish the entropy mixing condition.  Thus in terms of the
coupling analysis our conditions for variance and entropy mixing
are essentially the same.

Finally, we mention that our results actually hold
(with suitable modifications) for a much wider class of spin
systems on trees than just the Ising model, including the
Potts model and models with hard constraints such as the 
zero-temperature antiferromagnetic Potts model (proper colorings)
and the hard-core lattice gas model (independent sets).  
We briefly outline some of these extensions at the end of the paper;
full details can be found in a companion paper~\cite{companion}.

The remainder of the paper is organized as follows.
In Section~\ref{sec:prelims} we give some basic definitions and notation.
Then in Section~\ref{sec:spat} we define the spatial mixing conditions and
relate them to the spectral gap and log-Sobolev constant.
The mixing conditions in the scenarios of
interest for the spectral gap and the log-Sobolev constant
are verified in Sections~\ref{sec:gap} and~\ref{sec:logsob} respectively.
Finally, in Section~\ref{sec:extensions} we mention some
extensions of our results to other models of interest. 
The proofs of some technical lemmas omitted from the main text
are collected in a supplement, Section~\ref{sec:supplement}.

%%%%%%%%%%%%%%%%%% ACKNOWLEDGMENTS %%%%%%%%%%%%%%%%%%
%\vskip-0.3in\hbox{}
\section*{Acknowledgments}
%\vskip-0.05in
F.~Martinelli would like to thank the Miller Institute, the Dept.\ of
Statistics and the Dept.\ of EECS of the University of
California at Berkeley for financial support and warm hospitality.
We also wish to thank E.~Mossel and Y.~Peres for very interesting
discussions about reconstruction on trees and related topics.

%%%%%%%%%%%%%%%%%% PRELIMINARIES %%%%%%%%%%%%%%%%%%
%\vskip-0.3in\hbox{}
\section{Preliminaries}
\label{sec:prelims}
%\vskip-0.3in\hbox{}
\subsection{Gibbs distributions on trees}
%\vskip-0.05in
For $b\ge 2$, let $\Tree^b$ denote the infinite, rooted $b$-ary tree
(in which every vertex has $b$ children).  We will be concerned 
with (complete) finite subtrees $T$ of~$\Tree^b$; if $T$ has depth~$m$
then it has $n=(b^{m+1}-1)/(b-1)$ vertices, and its 
{\it boundary\/}~$\partial T$ consists of the children (in~$\Tree^b$)
of its leaves, i.e., $|\partial T|=b^{m+1}$.  We identify subgraphs
of~$T$ with their vertex sets, and write $E(A)$ for the edges
within a subset~$A$, and $\partial A$ for the boundary of~$A$
(i.e., the neighbors of~$A$ in $(T\cup \partial T)\setminus A$).

Fix an Ising spin configuration~$\tau$ on the infinite tree~$\Tree^b$.  
We denote by $\Omega_T^\tau$ the set of (finite) spin configurations 
$\sigma\in\{\pm 1\}^{T\cup\partial T}$ that agree with~$\tau$
on~$\partial T$; thus $\tau$ specifies a {\it boundary condition\/}
on~$T$.  Usually we abbreviate~$\Omega_T^\tau$ to~$\Omega$.
For any $\eta\in\Omega$ and any subset $A\subseteq T$, we denote
by~$\mu_A^\eta$ the Gibbs distribution over~$\Omega$ conditioned
on the configuration outside~$A$ being~$\eta$: i.e., if 
$\sigma\in\Omega$ agrees with~$\eta$ outside~$A$ then $$
   \mu_A^\eta(\sigma) \propto \exp\Bigl[\beta\bigl(\sum\nolimits_{xy\in E(A\cup\partial A)}\sigma_x\sigma_y + h\sum\nolimits_{x\in A}\sigma_x\bigr)\Bigr],  $$
where $\beta$ is the inverse temperature and $h$ the external field.
We define $\mu_A^\eta(\sigma)=0$ otherwise.  In particular, 
when $A=T$, $\mu_T^\tau$ is simply the Gibbs distribution
on the whole of~$T$ with boundary condition~$\tau$; we abbreviate
$\mu_T^\tau$ to~$\mu$.

For a function $f:\Omega\to\Rset$ we denote by 
$\mu_A^\eta(f)= \sum_{\sigma\in\Omega} \mu_A^\eta(\sigma)f(\sigma)$
the {\it expectation\/} of~$f$ w.r.t.\ the distribution~$\mu_A^\eta$.
It will be convenient to view $\mu_A^\eta(f)$ as a function
of~$\eta$, defined by $\mu_A(f)(\eta)=\mu_A^\eta(f)$, the
{\it conditional expectation\/} of~$f$.  Note that $\mu_A(f)$ is
a function from~$\Omega$ to~$\Rset$ but depends only on the 
configuration outside~$A$.
We write $\VarAeta(f)=\muAeta(f^2)-\muAeta(f)^2$ and (for $f\ge 0$)
$\EntAeta(f)=\muAeta(f\log f)-\muAeta(f)\log\muAeta(f)$ for
the {\it variance\/} and {\it entropy\/} of~$f$ respectively
w.r.t.\ $\muAeta$.  Note that $\VarAeta(f)=0$
iff, conditioned on the configuration outside~$A$ being~$\eta$, 
$f$~does not depend on the configuration inside~$A$.  The same
holds for $\EntAeta(f)$.  In case $A=T$ we use the abbreviations
$\mu(f), \Var(f)$ and $\Ent(f)$.

We record here some basic properties of variance and entropy
that we use throughout the paper:\par\noindent
{\bf (i)} For $B\subseteq A\subseteq T$,
\begin{equation}
  \label{e:var_dec}
  \Var^{\eta}_{A}(f) = \E^{\eta}_{A}[\Var_{B}(f)] +
   \Var^{\eta}_{A}[\E_{B}(f)].
\end{equation}
This equation expresses a decomposition of the variance into
the local conditional variance in~$B$ and the variance of the
projection outside~$B$.\par\noindent
{\bf (ii)} If $A=\bigcup_i A_i$ for disjoint~$A_i$, and the
Gibbs distribution~$\mu^{\eta}_{A}$ is the product of its marginals 
over the~$A_i$, then for any function~$f$,
\begin{equation}\label{eq:var_prod}
\Var^{\eta}_{A}(f) \le \sum_i \E^{\eta}_{A}[\Var_{A_i}(f)].
\end{equation} 
\par\noindent
{\bf (iii)} For any two subsets $A,B\subseteq T$ such that
$(\partial A) \cap B = \emptyset$, and for any function~$f$,
\begin{equation}
\label{e:convex}
\E[\Var_{A}(\E_{B}(f))] \le
\E[\Var_{A}(\E_{A\cap B}(f))] .
\end{equation}
Properties~(ii) and~(iii) are consequences of the fact that 
variance w.r.t.\ a fixed measure is a convex functional.
\par\smallskip\noindent
All three properties (i), (ii) and~(iii) also hold with $\Var$ replaced
by~$\Ent$.

%\vskip-0.3in\hbox{}
\subsection{The Glauber dynamics}
%\vskip-0.05in
The {\it Glauber dynamics\/} on $T$ with boundary
conditions $\t$ is the 
continuous time Markov chain on~$\Omega=\Omega_T^\tau$ with Markov
generator $\cL\equiv \cL_T^\t$ given by
\begin{equation}
  \label{generator}
  (\cL f)(\s) = \sum_{x\in T}c_x(\s) [f(\s^{x})-f(\s)], 
\end{equation}
where $\s^x$ denotes the configuration obtained from $\s$ by flipping
the spin at the site~$x$, and $c_x(\s)$ denotes the flip rate at~$x$.
Although all our results apply to any choice of finite--range, 
uniformly positive and bounded flip rates satisfying
the detailed balance condition w.r.t.\ the Gibbs
measure, for simplicity in the sequel we will work with a specific 
choice known as the {\it heat-bath\/} dynamics:
$$
c_x(\s)= \mu^{\s}_{\{x\}}(\s^x) = \frac{1}{1+w_x(\sigma)}, \quad \text{where}\quad
w_x(\sigma)=\exp\bigl[2\beta\sigma_x(\sum_{xy\in E}\sigma_y + h)\bigr].
$$

It is a well-known fact (and easily checked) that the Glauber 
dynamics is ergodic and reversible w.r.t.\ the Gibbs 
distribution~$\mu=\mu_T^\tau$, and so converges to the stationary
distribution~$\mu$.  The rate of convergence is often measured
using two concepts from functional analysis:
the {\it spectral gap\/} and the {\it logarithmic Sobolev constant}.
For a function
$f:\Omega\to\Rset$, define the {\it Dirichlet form\/} of~$f$ associated
with the generator~$\cL$ by
\begin{equation}\label{eq:df}
\df(f):= \half\sum_x \mu\bigl(c_x \bigl[f(\s^x)-f(\s)\bigr]^2\bigr)
   = \sum_x\mu(\Var_{\{x\}}(f)).
\end{equation}
(The l.h.s.\ here is the general definition for any choice of the
flip rates~$c_x$;
the last equality holds when specializing to the case of the heat-bath
dynamics.)
The \emph{spectral gap} $\cgap(\mu)$ and the \emph{logarithmic Sobolev
constant} $\csob(\mu)$ of the chain are then
defined by 
\begin{eqnarray}
   \cgap(\mu) = \inf_{f} {{\df(f)}\over\Var(f)};\qquad
   \csob(\mu) = \inf_{f\ge 0} {{\df(\sqrt{f}\,)}\over\Ent(f)},\label{eq:cgapcsob}
\end{eqnarray}
where the infimum in each case is over non-constant functions~$f$.

As is well known, these two quantities measure the rate
of exponential decay as $t\to \infty$ of the variance and relative
entropy respectively (see, e.g.,~\cite{Saloff}).
The quantity~$\cgap$ also has a natural interpretation as
the smallest positive eigenvalue of $-\cL$.

We make the following important note. When discussing the asymptotics of~$\csob$
(or~$\cgap$) for a fixed boundary condition~$\tau$, we think of the infinite
sequence of Gibbs distributions~$\set{\mu^\tau_T}$, where~$T$ ranges over all
finite complete subtrees of~$\Tree^b$. In particular, when we say that
$\csob(\mu)=\csob(\mu^\tau_T)=\Omega(1)$ we mean that there exists a finite
constant~$C>0$ such that for every~$T$ (or equivalently, for every $\mu\in
\set{\mu^\tau_T}$), $\csob(\mu)\ge 1/C$.

We close this section by recalling some well-known relationships
between the above constants and certain notions of mixing time
of the Glauber dynamics.   Define
$h_t^\s(\h) = \frac{P_t(\s,\h)}{\mu(\h)}$, where
\hbox{$P_t(\s,\h):=\nep{t\cL}(\s,\h)$} is the transition kernel at
time~$t$.  Then, for $1\le p \le \infty$,  define
\begin{equation}
  \label{T_p}
T_p := \min\Bigl\{t>0:\, \sup_\s \|h_t^\s-1\|_p \le \frac{1}{\nep{}}\Bigr\}  
\end{equation}
where $\|f\|_p$ denotes the $L^p(\O,\mu)$ norm of $f$. The time 
$T_1$ is usually called simply the \emph{mixing time} of the chain.
Standard results relating~$T_p$ to the spectral gap and log-Sobolev
constant (see, e.g.,~\cite{Saloff}), when specialized to the
Glauber dynamics, yield the following:
\begin{theorem}\label{th:backgnd}
On an $n$-vertex $b$-ary tree~$T$ with boundary condition~$\tau$,
\begin{itemize}
\item[{\rm(i)}] $\cgap(\mu)^{-1}\le T_1 \le \cgap(\mu)^{-1}\times C_1 n$;
\item[{\rm(ii)}] $\cgap(\mu)^{-1}\le T_2 \le \csob(\mu)^{-1}\times C_2 \log n$,
\end{itemize}
where $\mu=\mu_T^\tau$ and
$C_1,C_2$ are constants depending only on $b,\beta$ and~$h$.\qquad\square
\end{theorem}
Finally, we note that our choice of the heat-bath dynamics is not
essential.  Since changing to any other reversible local update rule 
(e.g., the Metropolis rule) affects $\csob$ and~$\cgap$ 
by at most a constant factor, our analysis applies to
any choice of Glauber dynamics.

%%%%%%%%%%%%%%%%%% SPATIAL MIXING %%%%%%%%%%%%%%%%%%
%\vskip-0.3in\hbox{}
\section{Spatial mixing conditions for spectral gap and log-Sobolev}
\label{sec:spat}
%\vskip-0.05in 
In this section we define a certain {\it spatial mixing condition\/} 
(i.e., a form of weak dependence between the spin at a site and
the configuration far from that site) for a Gibbs distribution~$\mu$, and prove
that this condition implies that $\cgap(\mu)=\Omega(1)$.
An analogous condition implies that $\csob(\mu)=\Omega(1)$.  
Our spatial mixing conditions have two main advantages
over those used previously: first, the conditions for the spectral gap and the
log-Sobolev constant are identical in form, allowing a uniform treatment;
second, and more importantly, 
they are {\it measure-specific}, i.e., they may hold for the Gibbs
distribution induced by some specific boundary configuration while not 
holding for other boundary configurations. Hence, the conditions are
sensitive enough to show rapid mixing for
specific boundaries even though the mixing time with other
boundaries is slow for the same choice of temperature and external
field. We also note that the results of this section hold not just for 
the Ising model but for any nearest-neighbor interaction model on a tree.

%\vskip-0.3in\hbox{}
\subsection{Reduction to block analysis}
%\vskip-0.05in
Before presenting the main result of this section, we need
some more definitions and background.  For each site~$x\in T$, 
let~$B_{x,\ell}\subseteq T$ denote the subtree (or ``block'') of
height~$\ell-1$ rooted at~$x$, i.e., $B_{x,\ell}$ consists of $\ell$~levels.
(If~$x$ is $k<\ell$ levels from the bottom of~$T$ then 
$B_{x,\ell}$ has only $k$~levels.) 
In what follows we will think of $\ell$ as a suitably large constant.
%and abbreviate~$B_{x,\ell}$ to~$B_x$.
By analogy with expression~(\ref{eq:df}) for the Dirichlet form,
let $\df_{\ell}(f)\equiv \sum_{x\in T} \E[\Var_{B_{x,\ell}}(f)]$ 
denote the local variation of~$f$ w.r.t.\ the blocks~$\set{B_{x,\ell}}$. 
A straightforward manipulation (see, e.g.,~\cite{Martinelli}, keeping
in mind that each site belongs to at most $\ell$ blocks) shows
that~$\cgap$ can be bounded as follows:
\begin{equation}
\label{e:block_to_single}
\cgap(\mu) \ge {1\over{\ell}}\cdot \inf_f \frac{\df_{\ell}(f)}{\Var(f)} \cdot
\min_{\eta,x} \cgap(\mu^{\eta}_{B_{x,\ell}}).
\end{equation}
As before, the infimum is taken over non-constant functions (and 
henceforth we omit explicit mention of this).  The importance
of~(\ref{e:block_to_single}) is 
that~$\min_{\eta,x}\cgap(\mu^{\eta}_{B_{x,\ell}})$
depends only on the size of~$B_{x,\ell}$ and~$\beta$, but not on the size of~$T$;
in fact, it is at least~$\Omega(e^{-c(b,\beta)\cdot\ell})$~\cite{BKMP}. 
Therefore, in order to show that~$\cgap$ is bounded by a constant 
independent of the size of~$T$, it
is enough to show that, for some finite~$\ell$,
$\Var(f)\le \const\times \df_{\ell}(f)$ for all functions~$f$.
This is what we will show below, under the relevant
spatial mixing condition.  
As a side remark, notice that
$\inf_f\frac{\df_{\ell}(f)}{\Var(f)}$ is exactly the spectral gap
of the Glauber dynamics based on flipping {\it blocks}~$B_{x,\ell}$,
rather than single sites~$x$.

An identical manipulation yields an analogous bound for
the log-Sobolev constant. 
For a non-negative function~$f$, let~$\dfe_{\ell}(f)\equiv \sum_{x\in
  T} \E[\Ent_{B_{x,\ell}}(f)]$. Then
\begin{equation}
\label{e:block_to_single_e}
\csob(\mu) \ge {1\over{\ell}}\cdot \inf_{f\ge 0} \frac{\dfe_{\ell}(f)}{\Ent(f)} \cdot
\min_{\eta,x} \csob(\mu^{\eta}_{B_{x,\ell}}).
\end{equation}
Hence to bound $\csob(\mu)$ it suffices to show that, for some
constant~$\ell$,
$\Ent(f)\le \const\times \dfe_{\ell}(f)$ for all~$f\ge 0$.

%\vskip-0.3in\hbox{}
\subsection{Spatial mixing}
%\vskip-0.1in
We are now ready to state our spatial mixing conditions, first for
the variance and then for the entropy.
For $x\in T$, write $T_x$ for the subtree rooted at~$x$, and
$\widetilde{T_x}$ for~$T_x\setminus{\set{x}}$, the subtree~$T_x$ 
excluding its root.
\begin{definition}
\label{d:vcond}
{\bf [Variance Mixing]}
We say that~$\mu=\mu^{\tau}_T$ satisfies $\vcond(\ell, \varepsilon)$ if for
every~$x\in T$, any~$\eta\in\Omega^{\tau}_{T}$ and any function~$f$ that does
not depend on~$B_{x,\ell}\,$, the following holds:
$$\Var^{\eta}_{T_x}[\E_{\widetilde{T_x}}(f)] \;\le\; \varepsilon \cdot
\Var^{\eta}_{T_x}(f). $$
\end{definition}

Let us briefly discuss the above condition. Essentially,
$\varepsilon=\varepsilon(\ell)$ gives the rate of decay with distance~$\ell$
of point-to-set correlations. To see this, note that the
l.h.s.~$\Var^{\eta}_{T_x}[\E_{\widetilde{T_x}}(f)]$ is the variance of 
the {\it projection\/} of~$f$ onto the root~$x$ of~$T_x$, which is at
distance~$\ell$ from the sites on which~$f$ depends. 
It is also worth noting that
the required uniformity in~$\eta$ in~$\vcond$ is not very restrictive: since the
distribution~$\mu^{\eta}_{T_x}$ depends only on the restriction of~$\eta$ to the
boundary of~$T_x$, and since~$\eta\in\Omega^{\tau}_T$ (i.e., $\eta$ agrees
with~$\tau$ on~$\partial T$ and therefore on the bottom boundary of~$T_x$), the
only freedom left in choosing~$\eta$ is in choosing the spin of the parent
of~$x$. Thus,~$\vcond$ is essentially a property of the distribution 
induced by
the boundary condition~$\tau$. It is this lack of uniformity (i.e., the fact
that we need not verify~$\vcond$ for other boundary conditions) that makes it
flexible enough for our applications. 

As the following theorem states,
if~$\vcond(\ell,\varepsilon)$ holds with~$\varepsilon\approx{1\over{2\ell}}$, 
then we get a lower bound on~$\cgap$:
\begin{theorem}
\label{t:spat_fast}
For any~$\ell$ and~$\delta>0$, if~$\mu$
satisfies~$\vcond(\ell,(1-\delta)/2(\ell+1-\delta))$ then
$\Var(f)\le{3\over\delta}\cdot\df_{\ell}(f)$ for all~$f$.  In particular,
%\marginpar{\tiny FABIO: Shouldn't we state the converses of both
%these theorems as separate theorems???}
if~$\vcond$ with the above parameters holds for some fixed~$\ell$ and
$\delta>0$, for all $\mu=\mu^\tau_T$ with~$T$ a full subtree, 
then $\cgap(\mu)=\Omega(1)$.  Conversely, if $\cgap(\mu)=\O(1)$ then for
all~$T$, $\mu^\tau_T$ satisfies $\vcond(\ell,c\nep{-\th \ell})$ for some
constants $c,\th>0$ and all~$\ell$.
\end{theorem}
\begin{remark*}
The second part of the theorem was already proved in~\cite{BKMP},
where it was shown that for general nearest-neighbor spin systems on
any bounded degree graph, if~$\cgap(\mu)$ is bounded independently
of~$n$ then~$\mu$ exhibits an exponential decay of point-to-set correlations
(i.e.,$~\vcond(\ell,c\exp(-\th
\ell))$ holds for all~$\ell$).  The authors of~\cite{BKMP} posed the
question of whether the converse is also true.
Theorem~\ref{t:spat_fast} (which holds for general nearest-neighbor
spin systems on a tree) answers this question affirmatively when the
graph is a tree. In fact, as is apparent from the above theorem, 
the decay of point-to-set correlations on a tree is either slower 
than linear or exponentially fast.
\end{remark*}

The analogous mixing condition for entropy and the log-Sobolev constant 
is the following:
\begin{definition}
\label{d:econd}
{\bf [Entropy Mixing]}
We say that~$\mu=\mu^{\tau}_T$ satisfies $\econd(\ell, \varepsilon)$ if for
every~$x\in T$, any~$\eta\in\Omega^{\tau}_T$ and any non-negative function~$f$
that does not depend on~$B_{x,\ell}\,$, the following holds:
$$\Ent^{\eta}_{T_x}[\E_{\widetilde{T_x}}(f)] \;\le\; \varepsilon \cdot
  \Ent^{\eta}_{T_x}(f). $$
\end{definition}

Before stating the analog of Theorem~\ref{t:spat_fast} relating~$\csob$ 
to~$\econd$, we need to define one more constant.
Let $p_{\mathrm{min}}=\min_{x,s,\eta\in\Omega^{\tau}_{T}}\mu^{\eta}_{T_x}(\sigma_x=s)$,
where~$s$ ranges over~$\set{+,-}$; i.e., $p_{\mathrm{min}}$
is the minimum probability of any spin value at any site with any
boundary condition.  It is easy to see that $p_{\mathrm{min}}\ge
{1\over 2}e^{-2\beta(b+|h|)}$, a constant depending only on $b,\beta,h$.
\begin{theorem}
\label{t:spat_fast_e}
For any~$\ell$ and~$\delta>0$, if~$\mu$
satisfies~$\econd(\ell,[(1-\delta)p_{\mathrm{min}}/(\ell+1-\delta)]^2)$ then
$\Ent(f)\le{2\over\delta}\cdot \dfe_{\ell}(f)$ for all $f\ge 0$.  In particular,
if~$\econd$ with the above parameters holds for some fixed~$\ell$ and
$\delta>0$, for all $\mu=\mu^\tau_T$ with~$\tau$ fixed and~$T$ an arbitrary full
subtree, then $\csob(\mu)=\O(1)$.  Conversely, if $\csob(\mu)=\O(1)$ then for
all~$T$, $\mu^\tau_T$ satisfies $\econd(\ell,c\nep{-\th \ell})$ for some
constants $c,\th>0$ and all~$\ell$.
\end{theorem}
In order to prove Theorems~\ref{t:spat_fast} and~\ref{t:spat_fast_e} 
it is convenient to work with
spatial mixing conditions that are somewhat more involved
than~$\vcond$ and~$\econd$. The main difference is that we want to allow for
functions that may depend on~$B_{x,\ell}$ (the first~$\ell$ levels of~$T_x$) and
thus need to introduce a term for this dependency. The modified conditions
express the property that the variance (entropy) of the projection of any
function~$f$ onto the root~$x$ of~$T_x$ can be bounded up to a constant factor
by the local variance (entropy) of~$f$ in~$B_{x,\ell}$, plus a negligible factor
times the local variance (entropy) of~$f$ in~$\widetilde{T_x}$.
As the following lemma states, the modified
conditions (with appropriate parameters) can be deduced from~$\vcond$
and~$\econd$.

\begin{lemma}\label{l:decay_vcond}
\begin{itemize}
\item[{\rm (i)}] For any~$\varepsilon<{1\over 2}$, if~$\mu=\mu^{\tau}_T$
  satisfies~$\vcond(\ell,\varepsilon)$ then for every~$x\in T$,
  any~$\eta\in\Omega^{\tau}_T$ and any function~$f$ we have
  $\Var^{\eta}_{T_x}[\E_{\widetilde{T_x}}(f)] \le {2-\varepsilon'\over
  1-\varepsilon'}\cdot\E^{\eta}_{T_x}[\Var_{B_{x,\ell}}(f)]  + {\varepsilon'\over
    1-\varepsilon'}\cdot\E^{\eta}_{T_x}[\Var_{\widetilde{T_x}}(f)]$,
  with~$\varepsilon'=2\varepsilon$.
\item[{\rm (ii)}] For any~$\varepsilon < p_{\mathrm{min}}^2$,
  if~$\mu=\mu^{\tau}_T$ satisfies~$\econd(\ell,\varepsilon)$ then for
  every~$x\in T$, any~$\eta\in\Omega^{\tau}_T$ and any function~$f\ge 0$ 
we have $\Ent^{\eta}_{T_x}[\E_{\widetilde{T_x}}(f)] \le {1\over
  1-\varepsilon'}\cdot\E^{\eta}_{T_x}[\Ent_{B_{x,\ell}}(f)] + {\varepsilon'\over
    1-\varepsilon'}\cdot\E^{\eta}_{T_x}[\Ent_{\widetilde{T_x}}(f)]$,
  with~$\varepsilon'={{\sqrt{\varepsilon}}\over p_{\mathrm{min}}}$.
\end{itemize}
\end{lemma}

\begin{remark*}
We note that with extra work, part~(ii) of Lemma~\ref{l:decay_vcond} can be
improved to hold with~$\varepsilon'=c(p_{\mathrm{min}})\varepsilon$. We give the
weaker bound because it is simpler to prove while still enough for our applications.
\end{remark*}
Similar statements to those in Lemma~\ref{l:decay_vcond} appeared
in~\cite{Cesi'}. We defer our proof to Section~\ref{sec:supplement}.

We can now prove Theorems~\ref{t:spat_fast} and~\ref{t:spat_fast_e} by working
with the modified spatial mixing conditions of Lemma~\ref{l:decay_vcond}.

\begin{proofof}{Theorems~\ref{t:spat_fast} and~\ref{t:spat_fast_e}}
  Here we only prove the forward direction of both theorems. The reverse
  direction of Theorem~\ref{t:spat_fast} was proved in \cite{BKMP}, as
  already mentioned above. The proof of the reverse direction of
  Theorem~\ref{t:spat_fast_e} is deferred to Section~\ref{sec:supplement}
  because it uses machinery developed later in the paper.
  
  The main step in the proof of the forward direction is to show the
  following claim:
  \begin{claim}\label{c:org_spat_fast}
  If for every~$x\in T$, any~$\eta\in\Omega^{\tau}_{T}$ and any function~$f$,
    $$\Var^{\eta}_{T_x}[\E_{\widetilde{T_x}}(f)] \;\le\; c\cdot
    \E^{\eta}_{T_x}[\Var_{B_{x,\ell}}(f)] \;+\; \left({1-\delta\over
        \ell}\right) \cdot \E^{\eta}_{T_x}[\Var_{\widetilde{T_x}}(f)],
    $$
    then~$\Var(f)\le{c\over\delta}\cdot\df_{\ell}(f)$ for all~$f$. The same
    implication holds when $\Var$ is replaced by~$\Ent$,~$\df_\ell$
    is replaced by~$\dfe_\ell$ and the function~$f$ is restricted to be
    non-negative.
  \end{claim}
  
  Observe that the hypothesis of
  Theorem~\ref{t:spat_fast} together with part~(i) of Lemma~\ref{l:decay_vcond}
  establishes the hypothesis of Claim~\ref{c:org_spat_fast} with~$c\le 3$, and
  similarly, the hypothesis of
  Theorem~\ref{t:spat_fast_e} together with part~(ii) of
  Lemma~\ref{l:decay_vcond} establishes the hypothesis of
  Claim~\ref{c:org_spat_fast} (after the necessary replacement
  of symbols) with~$c \le 2$.
  
  It therefore suffices to prove Claim~\ref{c:org_spat_fast}. We prove
  only the formulation with~$\Var$ and~$\df_\ell$ since
  the proof for the formulation with~$\Ent$ and~$\dfe_\ell$ is identical once we
  make the same replacements in the text of the proof.  As will be
clear below, the proof uses only  properties which are common to 
both~$\Var$ and~$\Ent$.
  
  Consider an arbitrary function~$f:\Omega\to\Rset$.  Our first goal is to
  relate $\Var(f)$ to the projections
  $\Var^{\eta}_{T_x}[\E_{\widetilde{T_x}}(f)]$ for $x\in T$, so that we can
  apply the spatial mixing condition of the hypothesis.  Recall that~$T$
  has~$m+1$ levels, and define the increasing sequence $\emptyset=F_0\subset
  F_1\subset\ldots\subset F_{m+1}=T$, where $F_i$ consists of all sites in the
  lowest $i$ levels of~$T$.  Thus $F_i$ is a forest of height~$i-1$.
  Using~(\ref{e:var_dec}) recursively, and the facts that
  $\mu_{F_{i+1}}(\mu_{F_{i}}(f)) = \mu_{F_{i+1}}(f)$ and $\mu_{F_0}(f)=f$, we
  obtain
%  \Var(f) = \sum\nolimits_{i=1}^{m+1}
%  \E[\Var_{F_i}(\E_{F_{i-1}}(f))]. $$
\begin{eqnarray*}
  \Var(f)&=& \E[\Var_{F_1}(f)] +
  \Var[\E_{F_1}(f)]\\
  &=&\E[\Var_{F_1}(f)] +
  \E[\Var_{F_2}(\E_{F_1}(f))] +
  \Var[\E_{F_2}(\E_{F_1}(f))]\\
  &\vdots&\\
  &=&\sum_{i=1}^{m+1}
  \E[\Var_{F_i}(\E_{F_{i-1}}(f))].
\end{eqnarray*}
Now a fundamental property of
nearest-neighbor interaction models on a tree is that, given the configuration
on~$T\setminus F_i$, the Gibbs distribution on~$F_i$ becomes a product of the
marginals on the subtrees rooted at the sites~$x\in F_i\setminus
F_{i-1}$. Using inequality~(\ref{eq:var_prod}) for the variance
of a product measure, we therefore have that
\begin{equation} \label{e:var_dec_proj}
 \Var(f) \le \sum_{i=1}^{m+1}
  \sum_{x\in F_i\setminus F_{i-1}}
  \E[\Var_{T_x}(\E_{F_{i-1}}(f))]
  \le  \sum_{x\in T}\E[\Var_{T_x}(\E_{\widetilde{T_x}} (f))],
\end{equation}
where in the second inequality we used the convexity of the
variance as in~(\ref{e:convex}).
    
Notice that so far we have not used the spatial mixing condition in the
hypothesis of Claim~\ref{c:org_spat_fast}, but only a natural martingale
structure induced by the tree.  Let us denote the final sum
in~(\ref{e:var_dec_proj}) by~$\PV(f)$. In order to bound~$\cgap$, we need to
compare the projection terms~$\Var_{T_x}(\E_{\widetilde{T_x}} (f))$ in~$\PV(f)$
with the local conditional variance terms in~$\df_{\ell}(f)$.  
For example, notice
that if~$\mu$ were the product of its single-site marginals
then~$\Var_{T_x}(\E_{\widetilde{T_x}} (f)) \le \E_{T_x}[\Var_x(f)]$
and~$\cgap=1$. However, in general the variance of the projection on~$x$ may
also involve terms which depend on other sites, and may lead to a factor that
grows with the size of~$T_x$.  We will use the spatial mixing condition in order
to preclude the latter possibility.  Specifically, we show that if for
every~$x\in T$, any~$\eta\in\Omega^{\tau}_{T}$ and any function~$g$,
$\Var^{\eta}_{T_x}[\E_{\widetilde{T_x}}(g)] \;\le\; c\cdot
\E^{\eta}_{T_x}[\Var_{B_{x,\ell}}(g)] + \varepsilon \cdot
\E^{\eta}_{T_x}[\Var_{\widetilde{T_x}}(g)]$ then for every~$x\in T$
and~$\eta\in\Omega$,
\begin{equation}
  \Var^{\eta}_{T_x}[\E_{\widetilde{T_x}} (f)] \;\le\; c\cdot \E^{\eta}_{T_x}[\Var_{B_{x}}(f)] 
  + \varepsilon\cdot\!\!\!\! \sum_{y \in B_{x}\cup\partialb B_{x},y\ne x}
%\mathop{\sum_{y \in B_{x}\cup\partialb B_{x}}}_{y\ne x}
  \E^{\eta}_{T_x}[\Var_{T_y}(\E_{\widetilde{T_y}} (f))],
  \label{e:block_bound}
\end{equation}
where we have abbreviated~$B_{x,\ell}$ to~$B_x$ and
$\partialb B_{x}$ stands for the boundary of~$B_{x}$ excluding the parent
of~$x$, i.e., the bottom boundary of~$B_{x}$.
Notice that the last term in~(\ref{e:block_bound}) is relevant 
only when~$x$ is at distance at
least~$\ell$ from the bottom of~$T$. When~$x$ belongs to one of
the~$\ell$ lowest levels of~$T$ then~$T_x=B_{x}$, and thus trivially
$\Var^{\eta}_{T_x}[\E_{\widetilde{T_x}} (f)] \le
\E^{\eta}_{T_x}[\Var_{B_{x}}(f)]$.

Let us assume~(\ref{e:block_bound}) for now and conclude the proof of 
the theorem.  Applying~(\ref{e:block_bound}) for every~$x$ and~$\eta$, 
and using the hypothesis that $\varepsilon={{1-\delta}\over{\ell}}$ and
the fact that each site appears in at most $\ell$~blocks, we get
\begin{eqnarray*}
  \PV(f) &\le & c\cdot \df_{\ell}(f) + \varepsilon \cdot \sum_{x\in T}
   \sum_{y \in B_{x}\cup\partialb B_{x},y\ne x}
%\mathop{\sum_{y \in B_{x}\cup\partialb B_{x}}}_{y\ne x}
  \E[\Var_{T_y}(\E_{\widetilde{T_y}} (f))]\\
  &\le& c\cdot \df_{\ell}(f) + \varepsilon \ell \cdot \sum_{y \in T}
  \E[\Var_{T_y}(\E_{\widetilde{T_y}} (f))]\\
  &=&   c\cdot \df_{\ell}(f) + (1-\delta) \PV(f),
\end{eqnarray*}
and hence
$$\Var(f) \;\le\; \PV(f) \;\le\; {c\over \delta}\cdot \df_{\ell}(f),$$
proving Claim~\ref{c:org_spat_fast}. 
We now return to proving~(\ref{e:block_bound}).

Let~$g=\E_{T_x\setminus (B_{x}\cup\partialb B_{x})}(f)$. Once we notice
that~$\E_{\widetilde{T_x}}(f) = \E_{\widetilde{T_x}}(g)$,
we can use the spatial mixing assumption that precedes~(\ref{e:block_bound}) to deduce
%  \Var^{\eta}_{T_x}[\E_{\widetilde{T_x}} (f)] \le c\cdot
%  \E^{\eta}_{T_x}[\Var_{B_{x}}(g)] +
%  \varepsilon \cdot \E^{\eta}_{T_x}[\Var_{\widetilde{T_x}}(g)]
%  \le  c\cdot \E^{\eta}_{T_x}[\Var_{B_{x}}(f)] + \varepsilon\cdot
%  \E^{\eta}_{T_x}[\Var_{\widetilde{T_x}}(g)],  $$
\begin{eqnarray*}
  \Var^{\eta}_{T_x}[\E_{\widetilde{T_x}} (f)] &\le& c\cdot
  \E^{\eta}_{T_x}[\Var_{B_{x}}(g)] +
  \varepsilon \cdot \E^{\eta}_{T_x}[\Var_{\widetilde{T_x}}(g)]\\
  &\le&  c\cdot \E^{\eta}_{T_x}[\Var_{B_{x}}(f)] + \varepsilon\cdot
  \E^{\eta}_{T_x}[\Var_{\widetilde{T_x}}(g)],
\end{eqnarray*}
where we used~(\ref{e:convex}) for the second inequality. We will be done once
we show that
\begin{equation}
  \E^{\eta}_{T_x}[\Var_{\widetilde{T_x}}(g)] \;\le\;
   \sum_{y \in B_{x}\cup\partialb B_{x},y\ne x}
%  \mathop{\sum_{y \in B_{x}\cup\partialb B_{x}}}_{y\ne x}
  \E^{\eta}_{T_x}[\Var_{T_y}(\E_{\widetilde{T_y}} (f))].
  \label{e:block_single_dec}
\end{equation}
But (\ref{e:block_single_dec}) follows from a similar argument to
that used earlier to show~$\Var(f)\le \PV(f)$, 
starting from the 
fact that $g=\mu_{F'_k}(f)$, where the forests $F'_i$ are defined
analogously to the~$F_i$ earlier but restricted to the subtree~$T_x$,
and $k=\height(x)-\ell$.  We omit the details.

This concludes the proof of Claim~\ref{c:org_spat_fast}, and thus of
Theorems~\ref{t:spat_fast} and~\ref{t:spat_fast_e}.
\end{proofof}

%%%%%%%%%%%%%%%%%% SPECTRAL GAP %%%%%%%%%%%%%%%%%%
%\vskip-0.3in\hbox{}
\section{Verifying spatial mixing for the spectral gap}\label{sec:gap}
%\vskip-0.05in
In this section, we will prove that the spectral gap of the Glauber dynamics is
bounded in all of the situations covered by Theorem \ref{A} in the
Introduction. 

In light of Theorem~\ref{t:spat_fast}, to bound the spectral gap it suffices to
verify the Variance Mixing condition~$\vcond(\ell,\varepsilon)$ with
$\varepsilon=(1-\delta)/2(\ell+1-\delta)$, for some constants $\ell,\delta>0$
independent of the size of~$T$.  In fact, we will show it with the asymptotically
tighter value $\varepsilon=c\exp(-\th \ell)$:
\begin{theorem}\label{t:var_decay}
In both of the following situations, there exists a positive constant~$\th$
(depending only on $b,\beta$ and~$h$) such that, for all~$T$,
the Gibbs distribution $\mu=\mu_T^\tau$ satisfies 
$\vcond(\ell, e^{-\vartheta\ell})$ for all~$\ell$\/{\rm :}
\begin{itemize}
\item[{\rm (i)}] $\tau$ is arbitrary, and either $\beta<\beta_1$ (with $h$
  arbitrary), or $|h|>h_c(\beta)$ (with $\beta$ arbitrary);
\item[{\rm (ii)}] $\tau$ is the \Plus-boundary condition, and
  $\beta,h$ are arbitrary.
\end{itemize}
As a corollary, in both situations~$\cgap(\mu)=\Omega(1)$.
\end{theorem}

\begin{remark*}
The validity of~$\vcond$, i.e, the decay of point-to-set correlations, is of
interest independently of its implication for the spectral gap (an implication
which is new to this paper): e.g., it is closely related to the purity of the
infinite volume Gibbs measure and to bit reconstruction problems on
trees~\cite{EKPS}.  In the special case of a free boundary and $h=0$,
part~{\rm (i)} of Theorem~\ref{t:var_decay} was first proved in~\cite{BRZ} via
a lengthy calculation, which was considerably simplified in~\cite{Ioffe}.  It
was later reproved in~\cite{BKMP} (for arbitrary boundary conditions) as a
consequence of the fact that the spectral gap is bounded in this situation.
An extension to general trees can be found in~\cite{EKPS} and~\cite{Ioffe2}.
Our motivation for presenting another proof of part~{\rm (i)} (in addition
to handling general fields~$h$) is the
simplicity of our argument compared with previous ones. As far as part~{\rm
  (ii)} is concerned, we are unaware of any 
previous results for the case of the
\Plus-boundary other than the fact that $\vcond(\ell,\varepsilon(\ell))$
must hold
with $\lim_{\ell\to \infty}\varepsilon(\ell)=0$ 
because the \Plus-phase is pure (see, e.g.,~\cite{Georgii}).
\end{remark*}

The rest of this section is divided into two parts.  First, we develop a 
general framework based on coupling in order to establish the exponential 
decay of point-to-set correlations.  This framework identifies two key
quantities, $\kappa$ and~$\gamma$, and states that when their product is
small enough then $\vcond$ holds.
Then, in the second part, we go back to proving
Theorem~\ref{t:var_decay} by calculating~$\kappa$ and~$\gamma$
for each of the above two regimes separately.

\subsection{A coupling argument for decay of point-to-set correlations}
\label{s:gap_coupling}
In this section we develop a coupling framework that enables us to verify the
exponential decay of point-to-set correlations from a simple calculation
involving single-spin distributions. 

First we need some additional notation.  When $x$ is not the root of~$T$, let
$\mu^+_{T_x}$ (respectively,~$\mu^{-}_{T_x}$) denote the Gibbs distribution in
which the parent of~$x$ has its spin fixed to~\Plus\/ (respectively,~\Minus\/)
and the configuration on the bottom boundary of~$T_x$ is specified by~$\tau$
(the global boundary condition on~$T$) \footnote{Notice that we do not specify
  the rest of the configuration outside~$T_x$ since it has no influence on the
  distribution inside~$T_x$ once the spin at the parent of~$x$ is fixed.
  However, since our distributions are defined over the whole configuration
  space, in the discussion below when the configuration outside~$T_x$ is
  relevant it will be understood from the context.}.  For two
distributions~$\mu_1$ and~$\mu_2$, we denote by~$\|\mu_1-\mu_2\|_x$ the
variation distance between the projections of~$\mu_1$ and~$\mu_2$ onto the spin
at~$x$. (Since the Ising model has only two spin values,
$\|\mu_1-\mu_2\|_x=|\mu_1(\sigma_x=+)-\mu_2(\sigma_x=+)|$.) Recall also that
$\eta^y$ denotes the configuration~$\eta$ with the spin at site~$y$ flipped.

We now identify two constants that are crucial for our coupling argument:
\begin{definition}\label{d:kp_gm}
For a sequence of Gibbs distributions~$\set{\mu^{\tau}_T}$ corresponding
to a fixed boundary condition~$\tau$, define $\kappa\equiv \kappa(\set{\mu^\tau_T})$
and $\gamma\equiv\gamma(\set{\mu^\tau_T})$ by
\begin{itemize}
\item[{\rm (i)}] $\kappa=\sup_T\max_z\|\mu^+_{T_z} -
  \mu^-_{T_z}\|_{z}$\/{\rm ;}
\item[{\rm (ii)}] $\gamma=\sup_T\max \|\mu^\eta_A -
  \mu^{\eta^y}_A\|_{z}$, where the maximum is taken over all
  subsets~$A\subseteq T$, all boundary configurations~$\eta$, all sites~$y$ on the
  boundary of~$A$ and all neighbors~$z\in A$ of~$y$.
\end{itemize}
\end{definition}
Note that $\kappa$ is the same as~$\gamma$, except that the maximization is
restricted to $A={{T_z}}$ and the boundary vertex~$y$ being the parent of~$z$;
hence always $\kappa\le\gamma$.  Since $\kappa$ involves Gibbs distributions
only on maximal subtrees~${{T_z}}$, it may depend on the boundary
condition~$\tau$ at the bottom of the tree.  By contrast, $\gamma$~bounds the
worst-case probability of disagreement for an {\it arbitrary\/} subset~$A$ and
arbitrary boundary configuration around~$A$, and hence depends only on $(\b,h)$
and not on~$\tau$.  It is the dependence of~$\kappa$ on~$\tau$ that opens up the
possibility of an analysis that is specific to the boundary condition. For
example, at very low temperature and with no external field,~$\kappa$ is close
to~$1$ in the free boundary case, while it is close to zero in the
\Plus-boundary case.

In our arguments~$\kappa$ will be used to bound the probability of a
disagreement percolating one level {\it down\/} the tree, namely, when we fix a
disagreement at~$x$ and couple the two resulting marginals on a child~$z$
of~$x$. On the other hand, $\gamma$~will be used in order to bound the
probability of a disagreement percolating one level {\it up\/} the tree, namely,
when we fix a single disagreement on the bottom boundary of a block, say at~$y$
(with the rest of the boundary configuration being arbitrary), and couple the
marginals on the parent of~$y$.

The novelty of our argument for establishing~$\vcond$ comes from the fact that
we identify {\it two separate\/} constants~$\kappa$ and~$\gamma$, and
consider their product, rather than working with~$\kappa$ alone:
\begin{theorem}\label{t:gap_coup}
  Any Gibbs distribution~$\mu=\mu^{\tau}_T$ satisfies~$\vcond(\ell,
  (\gamma\kappa b)^{\ell})$ for all~$\ell$, where~$\kappa$ and~$\gamma$ are the
  constants associated with the sequence~$\set{\mu^\tau_T}$ as specified in
  Definition~\ref{d:kp_gm}. In particular, if $\gamma\kappa b< 1$ then there
  exists a constant~$\th>0$ such that, for every~$T$, the measure
  $\mu=\mu^\tau_T$ satisfies $\vcond(\ell, e^{-\th\ell})$ for all $\ell$, and
  hence $\cgap(\mu)=\Omega(1)$.
\end{theorem}
\begin{proof}
  Fix arbitrary~$T$, $x\in T$, $\eta\in\Omega^{\tau}_T$. We need to show that
  for every function~$f$ that does not depend on~$B_{x,\ell}$,
  $\Var^{\eta}_{T_x}[\E_{\widetilde{T_x}}(f)] \le
  \varepsilon\cdot\Var^{\eta}_{T_x}(f)$ with~$\varepsilon=(\kappa\gamma
  b)^{\ell}$, i.e., projecting~$f$ onto the root (of~$T_x$)
  causes the variance to shrink by a factor~$\varepsilon$. As is well known,
  it is enough to establish a dual contraction, i.e., to
  consider an arbitrary function that depends only on the spin at the root and
  show that, when projecting onto levels~$\ell$ and below, the variance shrinks
  by a factor~$\varepsilon$. Formally, it is enough to show that for every
  function~$g$ that does not depend on~$\widetilde{T_x}$\footnote{Effectively
    this means that, conditioned on the configuration outside~$T_x$
    being~$\eta$,~$g$ depends only on the spin at the root~$x$.}  we have
  \begin{equation}\label{e:gap_coup}
    \Var^{\eta}_{T_x}[\E_{B_{x,\ell}}(g)] \;\le\; \varepsilon\cdot \Var^{\eta}_{T_x}(g).
  \end{equation}
  This is because for a function~$f$ that does not depend on~$B_{x,\ell}$, the
  variance of the projection can be written as
  $$\Var^{\eta}_{T_x}[\E_{\widetilde{T_x}}(f)] \;=\;
  \Cov^{\eta}_{T_x}(f,\E_{\widetilde{T_x}}(f)) \;=\;
  \Cov^{\eta}_{T_x}(f,\E_{B_{x,\ell}}(\E_{\widetilde{T_x}}(f))) \;\le$$
  $$\sqrt{\Var^{\eta}_{T_x}(f)\cdot
    \Var^{\eta}_{T_x}[\E_{B_{x,\ell}}(\E_{\widetilde{T_x}}(f)]} \enspace,$$
  where~$\Cov^\eta_A(f,f')$ denotes the
  covariance~$\E^\eta_A(ff')-\E^\eta_A(f)\E^\eta_A(f')$ and the last inequality is
  an application of Cauchy-Schwartz. We then have
  $$\Var^{\eta}_{T_x}[\E_{\widetilde{T_x}}(f)] \;\le\; \Var^{\eta}_{T_x}(f)
  \cdot{\Var^{\eta}_{T_x}[\E_{B_{x,\ell}}(\E_{\widetilde{T_x}}(f))]\over
    \Var^{\eta}_{T_x}[\E_{\widetilde{T_x}}(f)]}\enspace .$$
  If we
  assume~(\ref{e:gap_coup}) then the expression on the r.h.s. is bounded by
  $\varepsilon \cdot\Var^{\eta}_{T_x}(f)$ since $g=\E_{\widetilde{T_x}}(f)$ does
  not depend on~$\widetilde{T_x}$.
  
  We therefore proceed with the proof of~(\ref{e:gap_coup}), which goes via a
  coupling argument.  A {\it coupling\/} of two distributions~$\mu_1,\mu_2$
  on~$\Omega$ is any joint distribution~$\nu$ on~$\Omega^2$ whose marginals
  are~$\mu_1$ and~$\mu_2$ respectively.  For two configurations
  $\sigma,\sigma'\in\Omega$, let~$|\sigma - \sigma'|_{x,\ell}$ denote the {\it
    Hamming distance\/} between the restrictions of~$\sigma$ and~$\sigma'$
  to~$\widetilde{\partial}B_{x,\ell}$, i.e., the number of sites at
  distance~$\ell$ below~$x$ at which~$\sigma$ and~$\sigma'$ differ. Notice that
  $|\sigma-\sigma'|_{x,\ell}$ can be at most~$b^{\ell}$, the number of sites on
  the~$\ell$th level below~$x$.  Let~$\mu^{+}_{\widetilde{T_x}}$
  (respectively,~$\mu^{-}_{\widetilde{T_x}}$) stand for the Gibbs distribution
  where the spin at~$x$ is set to~$(+)$ (respectively,~$(-)$) and, as usual, the
  configuration on the bottom boundary of~$\widetilde{T_x}$ is specified
  by~$\tau$. Our goal will be to construct a coupling~$\nu$
  of~$\mu^+_{\widetilde{T_x}}$ and~$\mu^-_{\widetilde{T_x}}$ for which the
  expectation $\Expect_\nu|\sigma-\sigma'|_{x,\ell}\equiv
  \sum_{\sigma,\sigma'}\nu(\sigma,\sigma')|\sigma-\sigma'|_{x,\ell}$ is
  only~$(\kappa b)^\ell$.
  \begin{claim}
    \label{cl:decay_coup}
    For every~$x\in T$ and all~$\ell$ the following hold\/{\rm :}
    \begin{itemize}
    \item[{\rm (i)}] There is a coupling~$\nu$
      of~$\mu^{+}_{\widetilde{T_x}}$ and~$\mu^{-}_{\widetilde{T_x}}$ for
      which~$\Expect_\nu|\sigma-\sigma'|_{x,\ell} \le (\kappa b)^{\ell}$.
    \item[{\rm (ii)}] For any~$\eta,\eta'\in\Omega$ that have the same spin value at
      the parent of~$x$,
      $\|\mu^{\eta}_{B_{x,\ell}}-\mu^{\eta'}_{B_{x,\ell}}\|_x \le
      \gamma^{\ell}\cdot|\eta-\eta'|_{x,\ell}$. 
    \end{itemize}
  \end{claim}

  Let us assume Claim~\ref{cl:decay_coup} for the moment and complete the proof
  of~(\ref{e:gap_coup}). Consider an arbitrary~$g$ that does not depend
  on~$\widetilde{T_x}$. Let~$p=\mu^{\eta}_{T_x}(\sigma_x = +)$
  and~$q=1-p=\mu^{\eta}_{T_x}(\sigma_x = -)$. We also write~$g^{+}$
  for~$g(\sigma)$, where~$\sigma$ is any configuration that agrees with~$\eta$
  outside~$T_x$ and such that~$\sigma_x=+$. (This is well defined since~$g$ does
  not depend on~$\widetilde{T_x}$). We define~$g^-$ similarly. Without loss of
  generality we may assume that in the coupling~$\nu$ from
  Claim~\ref{cl:decay_coup} both the coupled configurations agree with~$\eta$
  outside~$T_x$ with probability~$1$. We then have
  \begin{eqnarray}
    \nonumber \Var^{\eta}_{T_x}[\E_{B_{x,\ell}}(g)] & = &\Cov^{\eta}_{T_x}[g,
    \E_{B_{x,\ell}}(g)]\\
    \nonumber & = &\Cov^{\eta}_{T_x}[g,
    \E_{\widetilde{T_x}}(\E_{B_{x,\ell}}(g))]\\
    \nonumber & = &pq(g^+-g^-)[\mu^{+}_{\widetilde{T_x}}(\E_{B_{x,\ell}}(g)) -
    \mu^{-}_{\widetilde{T_x}}(\E_{B_{x,\ell}}(g))]\\
    \nonumber & = & pq(g^+-g^-)\sum_{\sigma,\sigma'}
    \nu(\sigma,\sigma')[\E^{\sigma}_{B_{x,\ell}}(g) -
    \E^{\sigma'}_{B_{x,\ell}}(g)]\\
    & \le & pq|g^+-g^-|\sum_{\sigma,\sigma'}\nu(\sigma,\sigma')
    \|\mu^{\sigma}_{B_{x,\ell}} - \mu^{\sigma'}_{B_{x,\ell}}\|_x \cdot
    |g^+-g^-|\label{e:tv_gen}\\
    \nonumber & \le & pq(g^+-g^-)^2\sum_{\sigma,\sigma'}\nu(\sigma,\sigma')
    |\sigma - \sigma'|_{x,\ell}\cdot\gamma^{\ell}\\
    \nonumber & = & \gamma^\ell \cdot \Var^{\eta}_{T_x}(g) \cdot
    \Expect_{\nu}|\sigma-\sigma'|_{x,\ell}\\
    \nonumber & \le & (\gamma\kappa b)^{\ell}\cdot\Var^{\eta}_{T_x}(g).
  \end{eqnarray}
  In the sixth line here we have used part~(ii) of Claim~\ref{cl:decay_coup},
  and in the last line we have used part~(i).  This completes the
  proof of~(\ref{e:gap_coup}), and hence of Theorem~\ref{t:gap_coup}.
  We thus go back and prove Claim~\ref{cl:decay_coup}.
  
  The proof of Claim~\ref{cl:decay_coup} makes use of a standard recursive
  coupling along paths in the tree (as in, e.g.,~\cite{BKMP}).
  We start with part~(i), i.e., constructing a
  coupling~$\nu$ of~$\mu^{+}_{\widetilde{T_x}}$ and~$\mu^{-}_{\widetilde{T_x}}$
  with the required properties.  Since the underlying graph is a tree, we can
  couple~$\mu^{+}_{\widetilde{T_x}}$ and~$\mu^{-}_{\widetilde{T_x}}$
  recursively.  This goes as follows. First, given the spin at~$x$ the measures
  on~$T_z$ (where~$z$ ranges over the children of~$x$) are all independent of
  each other, so we can couple the projections on the~$T_z$'s independently.
  Then, we couple the two projections on~$T_z$ by first coupling the spin at~$z$
  using the optimal coupling (the one that achieves the variation distance) of
  the marginal measures on the spin at~$z$. Thus, the spins at~$z$ disagree with
  probability at most~$\kappa$. Once a coupled pair of spins at~$z$ is chosen,
  we continue as follows: if the spins at~$z$ agree then we can make the
  configurations in~$\widetilde{T_z}$ equal with probability~$1$ (because the
  two boundary conditions are the same); if the spins at~$z$ differ (i.e., one
  is~\Plus\ and the other~\Minus) then we recursively
  couple~$\mu^{+}_{\widetilde{T_z}}$ and $\mu^{-}_{\widetilde{T_z}}$.  We
  let~$\nu$ be the resulting coupling of $\mu^{+}_{\widetilde{T_x}}$ and
  $\mu^{-}_{\widetilde{T_x}}$, and notice that
  $\Expect_\nu|\sigma-\sigma'|_{x,l} \le (\kappa b)^{\ell}$ since for every
  site~$y$ at distance~$\ell$ below~$x$ the probability that the two coupled
  spins at~$y$ disagree is at most~$\kappa^{\ell}$.
  
  We go on to prove part~(ii) of Claim~\ref{cl:decay_coup}. First, by writing a
  telescopic sum and applying the triangle inequality we get that
  $$\|\mu^{\eta}_{B_{x,\ell}} - \mu^{\eta'}_{B_{x,\ell}}\|_x \;\le\;
  \sum_{i=1}^k \|\mu^{\eta^{(i-1)}}_{B_{x,\ell}} -
  \mu^{\eta^{(i)}}_{B_{x,\ell}}\|_x \enspace,$$
  where~$k=|\eta-\eta'|_{x,\ell}$
  and the sequence of configurations~$\eta^{(i)}$ is a site-by-site
  interpolation of the differences between~$\eta$ and~$\eta'$ in~$\partialb
  B_{x,\ell}$. (It suffices to interpolate only over the differences
  in~$\partialb B_{x,\ell}$ since the measure~$\mu^\eta_{B_{x,\ell}}$ depends
  only on the configuration in~$\partial B_{x,\ell}$ and since~$\eta$
  and~$\eta'$ agree on the parent of~$x$.)
%It is now enough to show that
%  $\|\mu^{\eta}_{B_{x,\ell}} - \mu^{\eta^y}_{B_{x,\ell}}\|_x \le \gamma^{\ell}$
%  for all~$\eta$ and~$y\in\partialb B_{x,\ell}$. This, however, follows by a
%  coupling argument as before, where this time we couple recursively along the
%  path from~$y$ to~$x$ (i.e., up the tree), and where the bound on the
%  probability of disagreement at each step is~$\gamma$.
It is now enough to show that
$\|\mu^{\eta}_{B_{x,\ell}} - \mu^{\eta^w}_{B_{x,\ell}}\|_x \le \gamma^{\ell}$
for all~$\eta$ and $w\in\partialb B_{x,\ell}$. This, however, follows by a
coupling argument as before, where this time we couple recursively
along the path from~$w$ to~$x$ (i.e., up the tree). Specifically, suppose by
induction that in our coupling there is already a path of disagreement going
from~$w$ to~$y$, where~$y$ is some site on the path from~$w$ to~$x$. Let~$z$
denote the parent of~$y$. At the next step we choose a coupled pair of spins
at~$z$ from the two distributions~$\mu^\eta_A$ and~$\mu^{\eta^y}_A$ (using an
optimal coupling for the projections onto the spin at~$z$), where the
subset~$A$ is~$B_{x,\ell}$ excluding the path from~$w$ to~$y$. The probability
of disagreement at~$z$ given the disagreement at~$y$ is then bounded
by~$\gamma$, by definition. If the resulting spins at~$z$ agree then the spins
on the rest of the path are coupled to agree with certainty, while if there 
is a disagreement at~$z$ we continue recursively starting from the disagreement
at~$z$. We therefore conclude that the probability of disagreement at~$x$ in
the resulting coupling is~$\gamma^\ell$, as required.
\end{proof}

\begin{remark*}
  We emphasize that Theorem~\ref{t:gap_coup} is not specific to the Ising
  model and generalizes to arbitrary nearest-neighbor models on a tree.
  Although we used the fact that the Ising model has only two possible spin
  values, the proof can easily be generalized to more than two spin values at the
  cost of a factor ${1\over p_{\mathrm{\min}}}$ in front of~$(\gamma\kappa
  b)^{\ell}$ in~$\vcond$, where~$p_{\mathrm{\min}}$ is the minimum probability
  of any spin value as defined just before Theorem~\ref{t:spat_fast_e}. Thus,
  since Theorem~\ref{t:spat_fast} also applies to general nearest-neighbor spin
  systems on a tree, we conclude that the implication from~$\gamma\kappa b < 1$
  to a bounded~$\cgap(\mu)$ holds for any such system (with the definitions
  of~$\kappa$ and~$\gamma$ extended in the obvious way to systems with more than
  two spin values).  
  The details can be found in the companion paper~\cite{companion}.
\end{remark*}

\subsection{Proof of Theorem~\ref{t:var_decay}}\label{s:gap_proof}
In this section we go back to proving Theorem~\ref{t:var_decay}. Using
Theorem~\ref{t:gap_coup}, all we need to do for the given choices of the Ising
model parameters is to bound~$\kappa$ and~$\gamma$ as in Definition~\ref{d:kp_gm}
such that~$\gamma\kappa b < 1$. In contrast to
Sections~\ref{sec:spat} and~\ref{s:gap_coupling}, which apply to general
nearest-neighbor spin systems on trees, here the calculations are specific to
the Ising model.

For both $\kappa$ and $\gamma$, we need to bound a quantity of the form
$\|\mu^{\eta}_A-\mu^{\eta^y}_A\|_{z}$, where $y\in\partial A$ and $z\in A$ is a
neighbor of~$y$.  The key observation is that this quantity can be expressed
very cleanly in terms of the ``magnetization'' at~$z$, i.e., the ratio of
probabilities of a \Minus-spin and a \Plus-spin at~$z$.  It will actually be
convenient to work with the magnetization {\it without\/} the influence of the
neighbor~$y$: thus we let $\mu^{\eta,y=*}_A$ denote the Gibbs distribution with
boundary condition~$\eta$, except that the spin at~$y$ is free (or equivalently,
the edge connecting~$z$ to~$y$ is erased). We then have:
\begin{proposition}\label{p:tv_mag}
  For any subset $A\subseteq T$, any boundary configuration~$\eta$, any site
  $y\in\partial A$ and any neighbor $z\in A$ of~$y$, we have
  $$\|\mu^{\eta}_A-\mu^{\eta^y}_A\|_{z} = K_{\beta}(R),$$
  where
  $R={\mu^{\eta,y=*}_A(\sigma_z=-)\over \mu^{\eta,y=*}_A(\sigma_z=+)}$ and the
  function $K_\beta$ is defined by
  $$K_{\beta}(a) = {1\over{e^{-2\beta}a+1}} - {1\over{e^{2\beta}a+1}}.  $$
%  $$K_{\beta}(a) = {e^{2\beta}a\over e^{2\beta}a+1} - {e^{-2\beta}a\over
%    e^{-2\beta}a+1}.  $$
\end{proposition}
\begin{proof}
  First, w.l.o.g.\ we may assume that the edge between~$y$ and~$z$ is the
  only one connecting~$y$ to~$A$; this is because a tree has no cycles, so once
  the spin at~$y$ is fixed $A$ decomposes into disjoint components that are
  independent.  We also assume w.l.o.g.\ that the spin at~$y$ is~\Plus\ 
  in~$\eta$, and we abbreviate~$\mu^\eta_A$ and~$\mu^{\eta^y}_A$ to~$\mu^{+}_A$
  and~$\mu^{-}_A$ respectively, and also~$\mu^{\eta,y=*}_A$ to~$\mu^{*}_A$.
  Thus $\|\mu^{\eta}_A-\mu^{\eta^y}_A\|_{z} = 
  |\mu^+_A(\sigma_z=+) - \mu^-_A(\sigma_z=+)|$, and
  $R={\mu^*_A(\s_z=-)\over \mu^*_A(\s_z=+)}$. 
  We write~$R^+$ for
  ${\mu^+_A(\s_z=-)\over \mu^+_A(\s_z=+)}$ and~$R^-$ for ${\mu^-_A(\s_z=-)\over
    \mu^-_A(\s_z=+)}$.  
  Since the only influence of~$y$ on~$A$ is through~$z$, we have
  $R^+=e^{-2\b}R$ and $R^-=e^{2\b}R$. The proposition now follows once we
  notice that, by definition of~$R^+$ and~$R^-$, $\mu^+_A(\sigma_z=+) =
  {1\over{R^+ + 1}}$ and $\mu^-_A(\sigma_z=+) = {1\over{R^- + 1}}$.
\end{proof}

Now it is easy to check that $K_\beta(a)$ is an increasing function in the
interval~$[0,1]$, decreasing in the interval~$[1,\infty]$, and is maximized at
$a=1$.  Therefore, we can always bound~$\kappa$ and~$\gamma$ from above by
$K_\beta(1)={e^{\beta}-e^{-\beta}\over e^{\beta}+e^{-\beta}}$. Indeed,
for~$\gamma$ we must make do with this crude bound because it has to 
hold for any
boundary configuration~$\eta$ and we cannot hope to gain by controlling the
magnetization~$R$. However, as we shall see,
for~$\kappa$ we can do better in some cases
by computing the magnetization at the root; when this differs from~$1$ 
we get a better bound than~$K_\beta(1)$.

We are now ready to proceed to the proof of Theorem~\ref{t:var_decay}:
\par\smallskip\noindent
{\bf(i) Arbitrary boundary conditions}
\nopagebreak\par\smallskip\noindent\nopagebreak
Here, the boundary condition~$\tau$ is arbitrary and we first consider
the (easy) case when~$\beta<\beta_0$ or~$|h|>h_c(\beta)$ (i.e.,~$h$ is
super-critical). In this case we do not need to resort to the calculation
of~$\kappa$ and~$\gamma$. As discussed in the Introduction, in this regime there
is a unique infinite volume Gibbs measure, so certainly the variation distance
at the root $\max_{\eta,\eta'}\|\mu^{\eta}_{B_{x,\ell}} -
\mu^{\eta'}_{B_{x,\ell}}\|_x$ goes to zero as~$\ell$ increases.  In fact, it is
not too difficult to see that in the above regime this variation distance goes to
zero exponentially fast, which directly implies the desired exponential decay of
correlations ($\vcond$) by plugging the bound on the variation distance into
expression~(\ref{e:tv_gen}) in the proof of Theorem~\ref{t:gap_coup}.

We go on to consider the more interesting regime when~$\beta_0\le \beta<\beta_1$
(i.e., intermediate temperatures) and the external field~$h$ is arbitrary.  Here
we use the fact that $\kappa \le \gamma \le K_\beta(1)$.
We then certainly have $\gamma\kappa b < 1$ whenever
$K_\beta(1)={e^{\beta}-e^{-\beta}\over e^{\beta}+e^{-\beta}}<\sqrt{1\over b}$,
i.e., whenever~$e^{-2\beta} > {\sqrt{b}-1\over \sqrt{b}+1}$.  From the
definition of~$\beta_1$ (see Section~\ref{sec:Ising_on_trees}),
this corresponds precisely to $\beta<\beta_1$.
(Observe how this non-trivial result drops out immediately from our
machinery, as expressed in the condition $\gamma^2 < {1\over b}$.)

This completes the verification of Theorem~\ref{t:var_decay} part~(i).
\par\smallskip\noindent
{\bf (ii) \Plus-boundary condition}
\par\smallskip\noindent
We now assume that~$\tau$ is the all-\Plus\ configuration and consider
arbitrary~$\beta$ and~$h$.  For convenience, we assume $h\ge -h_c(\beta)$ since
the case $|h|>h_c(\beta)$ was covered in part~(i) for all boundary
conditions~$\tau$. The important property of the regime $h\ge -h_c(\beta)$ is
that, for the \Plus-boundary, the spin at the root is at least as
likely to be~\Plus\ as it is to be~\Minus. We will show that $\gamma\kappa b < 1$
throughout this regime. Recall that we already showed that $\gamma \le K_\b(1) <
1$ for all finite~$\b$.  It is therefore enough to show that 
$\kappa\le {1\over b}$.

To calculate~$\kappa$, we need to bound the variation
distance $\|\mu^+_{T_z} - \mu^-_{T_z}\|_{z}$, which by
Proposition~\ref{p:tv_mag} is equal to~$K_{\beta}(R_z)$, where
$R_z={\mu^*_{T_z}(\s_z=-)\over \mu^*_{T_z}(\s_z=+)}$ and~$\mu^*_{T_z}$ is the
Gibbs distribution over the subtree~$T_z$ when it is disconnected from the rest
of~$T$ and the spins on its bottom boundary agree with~$\tau$. We thus have
$\kappa=\sup_T\max_{z\in T} K_{\beta}(R_z)$.

The final ingredient we need is a recursive computation of the
magnetization~$R_z$, the details of which (up to change of variables) can be
found in~\cite{Baxter} or~\cite{BRSSZ}. Let $y\prec z$ denote that~$y$ is a
child of~$z$. A simple direct calculation gives that $R_z=e^{-2\beta
  h}\prod_{y\prec z}F(R_y)$, where $F(a)\equiv F_{\beta}(a)={a+e^{-2\beta}\over
  e^{-2\beta}a+1}$.  In particular, if~$z$ is any site on the bottom-most level
of~$T$, then since the spins of the children of~$z$ are all set
deterministically to~\Plus, we get that $R_z=e^{-2\beta h}[F(0)]^b$.  We
thus define
\begin{equation}\label{e:J_def}
J(a) \equiv J_{\beta,h}(a) = e^{-2\beta h} [F(a)]^b
\end{equation}
and observe that, for any~$z\in T$, $R_z=J^{(\ell)}(0)$, where~$J^{(\ell)}$
stands for the~$\ell$-fold composition of~$J$, and $\ell$ is the
distance of~$z$ from the bottom boundary of~$T$.

\begin{figure}[h]
\centerline{\psfig{file=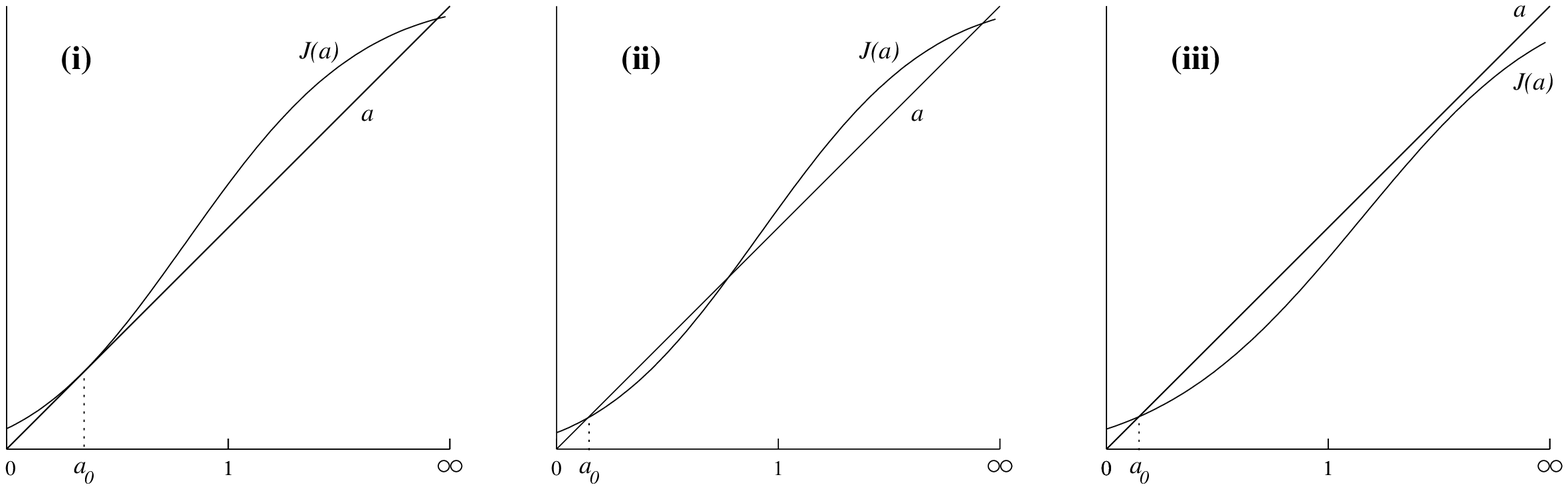,height=2in}}
\caption{Curve of the function $J(a)$, used in the proof of
Theorem~\ref{t:var_decay},
for $\beta>\beta_0$ and various values of the external field~$h$.
(i)~$h=-h_c(\beta)$; (ii)~$h_c(\beta)>h>-h_c(\beta)$; (iii)~$h>h_c(\beta)$.
The point~$a_0$ is the smallest fixed point of~$J$.}
\label{fig:J1}
\end{figure}

We now describe some properties of~$J$ that we use (refer to Fig.~\ref{fig:J1}):
$J$ is continuous and increasing on~$[0,\infty)$, with~$J(0)=e^{-2\beta(h+b)}>0$
and~$\sup_a J(a)=e^{-2\beta(h-b)}<\infty$. This immediately implies that $J$ has
at least one fixed point in $[0,\infty)$; we denote by~$a_0$ the least fixed
point.  Since~$a_0$ is the least fixed point and $J(0)>0$ then clearly $J'(a_0)
\le 1$, where $J'(a)\equiv {\partial J(a)\over \partial a}$ is the derivative
of~$J$. We also note that $a_0\le 1$ when $h\ge -h_c(\b)$, which corresponds to
the fact that for the \Plus-boundary and the above regime of~$h$, the spin
at the root is at least as likely to be \Plus\ as \Minus.

Now, since~$J$ is monotonically increasing and $a_0$ is the least fixed point
of~$J$, clearly $J^{(\ell)}(0)$ converges to~$a_0$ from below, i.e.,
$R_z\le a_0$ for every $z\in T$. Thus, since $a_0\le 1$ for $h\ge -h_c(\beta)$,
and the function~$K_{\beta}(a)$ is monotonically increasing in the
interval $[0,1]$, $K_{\beta}(R_z)\le K_{\beta}(a_0)$ for every $z\in T$.

What remains to be shown is that~$K_{\beta}(a_0) \le {1\over b}$. 
This follows from the fact that $J'(a_0)\le 1$,
together with the following lemma:
\begin{lemma}
\label{l:deriv_coup}
Let~$a_0$ be any fixed point of~$J$. 
Then~$K_{\beta}(a_0) = {1\over b}\cdot J'(a_0)$.
\end{lemma}
\begin{proofnobox}
{}From the definitions of~$J$ and~$F$ we have:  $$
\eqalign{
J'(a_0) &= e^{-2\beta h}\cdot b\cdot [F(a_0)]^{b-1}F'(a_0)\cr
& =  b\cdot J(a_0) \cdot {F'(a_0) \over F(a_0)}\cr
& =  b\cdot a_0 \cdot {F'(a_0) \over F(a_0)}\cr
& =  b\cdot a_0 \cdot \Bigl[{{1-e^{-4\beta}}\over{(a_0+e^{-2\beta})(e^{-2\beta}a_0+1)}}\Bigr]\cr
& =  b\cdot K_{\beta}(a_0).\qquad\square  \cr}  $$
\end{proofnobox}
This completes the verification of Theorem~\ref{t:var_decay} part~(ii).

%%%%%%%%%%%%%%%%%% LOG-SOBOLEV %%%%%%%%%%%%%%%%%%
%\vskip-0.3in\hbox{}
\section{Verifying spatial mixing for log-Sobolev}
\label{sec:logsob}
%\vskip-0.05in
In this section we will prove a uniform lower bound (independent of~$n$)
on the logarithmic Sobolev constant $c_{\rm sob}(\mu)$ in all the situations
covered by Theorem~\ref{B} in the Introduction.  

In light of Theorem~\ref{t:spat_fast_e}, to show $\csob=\Omega(1)$ we need only
prove the validity of the Entropy Mixing condition $\econd (\ell,
[(1-\delta)p_{\mathrm{min}}/2(\ell+1-\delta)]^2)$ for some constants~$\ell$ and
$\delta$ independent of the size of~$T$. In order to establish~$\econd$ in the
situations covered by Theorem~\ref{B}, we extend the coupling framework
developed in Section~\ref{s:gap_coupling} so that it can be used to
establish~$\econd$.  As before, we will use a condition on the
constants~$\kappa$ and~$\gamma$, which were defined in
Section~\ref{s:gap_coupling}. In fact, the condition on~$\kappa$ and~$\gamma$
for establishing~$\econd$ is practically the same as the one that was used to
establish~$\vcond$, which immediately transfers our~$\Omega(1)$ bound on~$\cgap$
for the relevant parameters to an~$\Omega(1)$ bound on~$\csob$ for the same
choice of parameters. The main result of this section is the following
relationship between $(\kappa,\gamma)$ and~$\econd$.

\begin{theorem}\label{t:em_coup}
  Any Gibbs distribution $\mu=\mu^{\tau}_T$ satisfies $\econd(\ell,
  c(\gamma\alpha)^{\ell/5})$ for all~$\ell$, where $\alpha=\max\set{\kappa b,
    1}$, $\kappa$ and~$\gamma$ are the constants associated with the
  sequence~$\set{\mu^\tau_T}$ as specified in Definition~\ref{d:kp_gm}, and~$c$
  is a constant that depends only on~$(b,\b,h)$. In particular, if
  $\max\set{\gamma\kappa b, \gamma} < 1$ then there exists a constant~$\th$ such
  that, for every~$T$, the measure $\mu=\mu^\tau_T$ satisfies $\econd(\ell,
  ce^{-\th\ell})$ for all $\ell$, and hence $\csob(\mu)=\Omega(1)$.
\end{theorem}
\begin{remark*}We should note that the above theorem, like its
counterpart for the spectral gap, holds for any spin system on a tree
(with the definitions of~$\kappa$ and~$\gamma$ generalized appropriately).
See the companion paper~\cite{companion} for details.
\end{remark*}

Since in Section~\ref{s:gap_proof} we have already calculated~$\kappa$
and~$\gamma$ for the regimes of interest and shown that in both cases
$\max\set{\gamma\kappa b,\gamma} < 1$, we have:
\begin{corollary}
\label{main_fabio}
In both of the following situations, 
%there exist positive constants $c$ and
%$\th$ (depending only on $b,\beta$ and~$h$) such that, for all~$T$, the Gibbs
%distribution $\mu\equiv\mu_T^\tau$ satisfies $\econd(\ell, ce^{-\vartheta\ell})$
%for all~$\ell$, and hence 
$\csob(\mu)=\Omega(1)$:
\begin{itemize}
\item[{\rm (i)}] $\tau$ is arbitrary, and either $\beta<\beta_1$ (with $h$
  arbitrary), or $|h|>h_c(\beta)$ (with $\beta$ arbitrary);
\item[{\rm (ii)}] $\tau$ is the \Plus-boundary condition and $\beta,h$
are arbitrary.
\end{itemize}
\end{corollary}
This completes the proof of our second main result, Theorem~\ref{B} stated
in the Introduction.

The first step in proving Theorem~\ref{t:em_coup} is a reduction of~$\econd$ to
a certain {\em strong concentration} property of~$\mu$, the Gibbs measure under
consideration. We believe that this concentration property, as well as
its connection to~$\econd$, may be of independent interest. The statement of this
property and the reduction of~$\econd$ to it is the content of
Section~\ref{subsec:econd_conc}. Then, in Section~\ref{subsec:em_kap_gam}, we complete the
proof of Theorem~\ref{t:em_coup} by relating the strong concentration property
to~$\kappa$ and~$\gamma$.

It is worth mentioning that we are also able to establish a general (but cruder)
bound on~$\csob$ as a function of~$\cgap$. Specifically, we can show that
$\csob = \Omega(1/\log n)\times \cgap$. Although we do not need this bound in this
paper, we present it in Section~\ref{subsec:crude} for future reference since
its proof is simple and short.

\subsection{Establishing $\econd$ via a strong concentration property.}
\label{subsec:econd_conc}
In this subsection we reduce~$\econd$ to a certain strong concentration property
of~$\mu$. In the next subsection, we will then establish this strong
concentration property as a function of~$\kappa$ and~$\gamma$ in order to prove 
Theorem~\ref{t:em_coup}. For simplicity and without loss of generality, we will
analyze the entropy mixing condition only for $T_x=T$ (the whole tree),
with root~$r$.

Let $\mu^+_{\widetilde T}$ and $\mu^-_{\widetilde T}$ denote the
Gibbs distributions on~${\widetilde T}$ with the spin at the root~$r$
set to~\Plus\ and~\Minus\ respectively (the boundary condition on
the leaves of~$T$ being specified by~$\tau$).  Define
\begin{equation*}
%\label{g}
g_+(\s) = \frac{\mu^+_{\widetilde T}(\s)}{\mu(\s)}= \left \{
\begin{array}{ll}
{1/ p} & \mbox{if $\sigma_r=(+)$,}\\
0 & \mbox{otherwise,}
\end{array}\right.
\end{equation*}
where $p=\mu(\s_r= +)$.
The key quantity we will work with in the sequel is the following:
\begin{equation*}
%\label{gl}
     g_+^{(\ell)}= \mu_{B_{r,\ell}}(g_+).
\end{equation*}
Note that $g_+^{(\ell)}(\s)$ depends only on the spins in~$\partial B_{r,\ell}$.
Indeed, let~$\sigma_{r,\ell}$ stand for the restriction of~$\sigma$ to~$\partial
B_{r,\ell}$, i.e., to the sites at distance~$\ell$ below~$r$.  It is easy to
verify that~$g_+^{(\ell)}(\sigma)$ is equal to
${\mu^{+}_{\widetilde{T}}(\sigma_{r,\ell})\over \mu(\sigma_{r,\ell})}$.  Thus,
for a given configuration~$\sigma$, $g_+^{(\ell)}(\sigma)$ is the ratio of the
probabilities of seeing the spins of~$\sigma$ at level~$\ell$ below the root~$r$
when the spin at~$r$ is~\Plus\ and when there is no condition on the spin
at~$r$, respectively. We define~$g_-$ and~$g_-^{(\ell)}$ in an analogous way.

The role played by the functions~$g_+^{(\ell)}$ and~$g_-^{(\ell)}$ is embodied
in the following theorem, which says that if these functions are sufficiently
tightly concentrated around their common mean value of~1 then the entropy mixing
condition~$\econd$ holds.
\begin{theorem}\label{th:conc_econd}
There exists a constant~$c$ (depending only on~$b$, $\b$ and~$h$) such that,
for any~$\delta\ge 0$, if
\begin{equation}\label{eq:strconc}
\mu\bigl[|g_s^{(\ell)} - 1| > \delta\bigr] \;\le\;
e^{-2/\delta}
\end{equation}
for $s\in\{+,-\}$, then we have $\Ent\bigl[\mu_{\widetilde T}(f)\bigr] \le
c \delta\Ent(f)$ for any non-negative function~$f$ that does not depend
on~$B_{r,\ell}\,$; in particular, $\econd(\ell,c\delta)$ holds.
\end{theorem}

\begin{proof}
Fix $\ell <m$ and a non-negative function $f$ that does not depend on
the spins inside the block~$B_{r,\ell}$. Since $\Ent(f')\le \Var(f')/\E(f')$
for every non-negative function~$f'$ (see, e.g.,~\cite{Saloff}) then
\begin{align}
  \Ent\bigl[\E_{\widetilde T}(f)\bigr] &\;\le \; \frac{\Var\bigl[\E_{\widetilde
      T}(f)\bigr]}{\E\bigl[\E_{\widetilde{T}}(f)\bigr]} \;=\; \ov {\E(f)}\cdot \Bigl[p
  \bigl[\E^+_{\widetilde T}(f) - \E(f)\bigr]^2 + (1-p) \bigl[\E^-_{\widetilde
      T}(f) - \E(f)\bigr]^2\Bigr]
  \nonumber \\
  &=\;\ov {\E(f)} \cdot \Bigl[p\Cov\bigl(g_+,f\bigr)^2 +
  (1-p)\Cov\bigl(g_-,f\bigr)^2\Bigr] \;\le\; \max_{s\in\set{+,-}}
  \frac{\Cov\bigl(g_s,f\bigr)^2}{\E(f)}\,, \label{th.1}
\end{align}
where~$\Cov$ denotes covariance w.r.t~$\mu$. Now observe that, since~$f$ does
not depend on~$B_{r,\ell}$, when computing
the covariance term in~(\ref{th.1}) the function~$g_s$ can be replaced
by~$g_s^{(\ell)}$, which depends only on the spins in~$\partial B_{r,\ell}$.
Thus, if we can show that~(\ref{eq:strconc}) implies
\begin{equation}\label{+}
\Cov\bigl(g_s^{(\ell)},f\bigr)^2  \le c\delta
\E(f)\Ent(f)
\end{equation}
for some constant~$c$, then by plugging~(\ref{+}) into~(\ref{th.1}) we will get
that $\Ent\bigl[\E_{\widetilde T}(f)\bigr] \le c\delta \Ent(f)$, as required.

To establish~(\ref{+}) we make use of the following technical lemma,
whose proof can be found in Section~\ref{sec:supplement}.

\begin{lemma}
\label{En_in_gen}
Let $\{\O,\cF,\nu\}$ be a probability space and let $f_1$ be a
mean-zero random variable such that $\ninf{f_1}\le 1$ and 
$\nu[\,|f_1|> \d\,]\le \nep{-2/\d}$ for some $\d\in (0,1)$. 
Let $f_2$ be a probability
density w.r.t.~$\nu$, \ie $f_2\ge 0$ and $\nu(f_2)=1$. 
Then there exists a numerical constant
$c'>0$ independent of $\nu,\, f_1,\,f_2$ and $\d$, such that
$\nu(f_1f_2)^2 \le c' \, \d \Ent_\nu(f_2)$. 
\end{lemma}
We apply this lemma with $\nu=\mu$ and
$$
f_1=\frac{\bigl(g_s^{(\ell)}-1\bigr)}{\ninf{g_s^{(\ell)}}};\quad
f_2=\frac{f}{\mu(f)}, $$
to deduce $\Cov\bigl(g_s^{(\ell)},f\bigr)^2 \le
c'\delta\ninf{g_s^{(\ell)}}^2 \mu(f)\Ent(f)$.  Noting also that
$\ninf{g_s^{(\ell)}} \le \ninf{g_s} \le 1/p_{\min}$, where~$p_{\min}$ was
defined just before Theorem~\ref{t:spat_fast_e}, this establishes~(\ref{+}) with
$c={c' / p_{\min}^2}$ and thus completes the proof of the theorem.
\end{proof}

\subsection{Proof of Theorem~\ref{t:em_coup}}\label{subsec:em_kap_gam}
In light of Theorem~\ref{th:conc_econd}, to prove Theorem~\ref{t:em_coup} it is
sufficient to verify the strong concentration property~(\ref{eq:strconc}) of the
functions~$g_s^{(\ell)}$ with $\delta=(\gamma\alpha)^{\ell/5}$.

In order to do this we appeal to a strong concentration of the Hamming distance
under the coupling~$\nu$ of~$\mu^+_{\widetilde{T}}$ and~$\mu^-_{\widetilde{T}}$,
as defined in the proof of Claim~\ref{cl:decay_coup}. Recall the notation used
in that claim, and notice that the Hamming distance is dominated by the size
of the population in the~$\ell$th generation of a specific branching process.
The following tail bound can be obtained using standard techniques from
the analysis of branching processes, and we defer the proof to the end
of this section.
\begin{lemma}\label{l:ham_conc}
  Let $\alpha=\max\set{\kappa b, 1}$. Then for every~$C>0$,
  $$
  \Pr_\nu \left[|\sigma-\sigma'|_{r,\ell} > C\alpha^\ell\right] \;\le \;
  e^{{1\over \ell+1}\bigl(1-{C\over 2e}\bigr)}.$$
\end{lemma}

\begin{corollary}\label{cr:g_coup_conc}
  For every~$C>0$ and $s\in\set{+,-}$,
  $$
  \Pr_\nu \left[\left|g_s^{(\ell)}(\sigma)-g_s^{(\ell)}(\sigma')\right| > C
    (\gamma\alpha)^\ell\right] \;\le \; e^{{1\over \ell+1}\bigl(1-{p_{\min} C\over
      2e}\bigr)}.$$
\end{corollary}
\begin{proof}
It is enough to show that
\begin{equation}\label{e:g_l_grad}
|g_s^{(\ell)}(\sigma)-g_s^{(\ell)}(\sigma')|\le {\gamma^{\ell} \over p_{\min}}
 \cdot |\sigma-\sigma'|_{r,\ell} 
\end{equation}
since we can then apply Lemma~\ref{l:ham_conc} with~$C$ replaced
by~$p_{\min}C$.  On the other hand,~(\ref{e:g_l_grad}) follows from part~(ii) of
Claim~\ref{cl:decay_coup} once we recall that
$g_s^{(\ell)}(\sigma)=\mu^\sigma_{B_{r,\ell}}(g_s)$ and that~$g_s$ depends only on the
spin at the root, implying that 
$|g_s^{(\ell)}(\sigma)-g_s^{(\ell)}(\sigma')|\le
\|\mu^\sigma_{B_{r,\ell}}-\mu^{\sigma'}_{B_{r,\ell}}\|_r
\cdot\ninf{g_s} \le \gamma^\ell \,|\sigma-\sigma'|_{r,\ell}\,/p_{\min}$.
\end{proof}

Before we go on with the proof of Theorem~\ref{t:em_coup}, let us compare
the way we used the constants~$\kappa$ and~$\gamma$ in the proof of
Corollary~\ref{cr:g_coup_conc} to the way we used them in the proof of
Theorem~\ref{t:gap_coup}. In both cases we used~$\kappa$ and~$\gamma$ to get
bounds for coupling ``down'' and ``up'' the tree respectively.  
Specifically, we
used~$\kappa$ to deduce that the Hamming distance between the coupled
configurations at the~$\ell$th level is about $(\kappa b)^\ell$, and we then
used~$\gamma$ to bound the effect of each discrepancy at the~$\ell$th level on
the spin at the root (or equivalently, on~$g_s^{(\ell)}$) by 
roughly~$\gamma^\ell$.  While in Theorem~\ref{t:gap_coup} it
was enough that the {\em average} Hamming distance when coupling down the tree
was bounded by $(\kappa b)^\ell$, here we need that this distance is not much
larger than $(\kappa b)^\ell$ {\em with high probability}.

We now return to the proof of Theorem~\ref{t:em_coup}. 
W.l.o.g.\ we may assume that
$\gamma\alpha \le 1$ since $\econd(\ell,1)$ always holds, and also that
$\gamma\alpha>0$ since if $\gamma=0$ then $\econd(\ell,0)$ holds
because then the spin
at the root~$r$ is independent of the rest of the configuration. Let
$a=(\gamma\alpha)^{-1}\ge 1$. Recall that we 
wish to establish~(\ref{eq:strconc}) with 
$\delta=a^{-\ell/5}$ for all large enough~$\ell$.
We will show only that
\begin{equation}\label{e:pos_conc}
\mu\bigl[g_s^{(\ell)}-1 > \delta\bigr] \le {1\over 2}e^{-2/ \delta}
\end{equation}
since the same bound on the negative tail can be achieved by an
analogous argument.

We start by applying Corollary~\ref{cr:g_coup_conc} with~$C=a^{\ell/4}$ to get
that, for every~$\varepsilon>0$,
\begin{equation}\label{e:conc_g_ham}
  \mu^s_{\widetilde{T}}\bigl[g_s^{(\ell)}-1 >
  \varepsilon\bigr] \;\le\;
\mu\bigl[g_s^{(\ell)}-1 >
  \varepsilon-a^{-3\ell/4}\bigr] + A,
\end{equation}
where $A=e^{{1\over \ell+1}\bigl(1-{p_{\min}a^{\ell/4}\over 2e}\bigr)}$ and we
have used the fact that~$\mu$ is a convex combination of~$\mu^+_{\widetilde{T}}$
and~$\mu^-_{\widetilde{T}}$.

Next, we notice that by definition of~$g_s^{(\ell)}$,
\begin{equation}\label{e:g_bnd}
  \mu^s_{\widetilde{T}}\bigl[g_s^{(\ell)}-1 >
    \varepsilon\bigr]  \; \ge \;
    (1+\varepsilon)\mu\bigl[g_s^{(\ell)}-1 >
      \varepsilon\bigr].
    \end{equation}
Combining~(\ref{e:conc_g_ham}) and~(\ref{e:g_bnd}) we get that, for every
$\varepsilon> 0$,
\begin{equation}\label{e:comb}
\mu\bigl[g_s^{(\ell)}-1 >
  \varepsilon\bigr] \;\le\; \left({1\over
    1+\varepsilon}\right)\left( \mu\bigl[g_s^{(\ell)}-1 >
    \varepsilon-a^{-3\ell/4}\bigr] + A\right).
\end{equation}
This immediately yields that, for every non-negative integer~$k$ and
$\varepsilon> 0$,
\begin{equation}\label{e:conc_interval}
  \mu\bigl[g_s^{(\ell)}-1>\varepsilon + ka^{-3\ell/4}\bigr]
  \;\le\; (1+\varepsilon)^{-(k+1)} + A\left({1+\varepsilon\over \varepsilon}\right),
\end{equation}
where we applied~(\ref{e:comb}) $k+1$ times, each time
increasing~$\varepsilon$ by~$a^{-3\ell/4}$.

Inequality~(\ref{e:pos_conc}) then follows (assuming~$\ell$ is large enough) by
applying~(\ref{e:conc_interval}) with $\varepsilon = a^{-\ell/4}$ and
$k=\lceil a^{\ell/2}\rceil$. This concludes the proof of
Theorem~\ref{t:em_coup}.
\stopproof\ignorespaces\bigbreak

Finally, we supply the missing proof of Lemma~\ref{l:ham_conc}.
\begin{proofof}{Lemma~\ref{l:ham_conc}}
  First notice that, by an exponential Markov inequality, it is enough to show
  that $\Expect_\nu \left[e^{t|\sigma-\sigma'|_{r,\ell}}\right] \le
  e^{2et\alpha^\ell}$ for all $t\le (2e(\ell+1)\alpha^\ell)^{-1}\le 1$.  We thus
  fix~$t$ as above and let~$D_{x,i}=\Expect_\nu
  \left[e^{t|\sigma-\sigma'|_{x,i}}\right]$, where~$\nu$ is the coupling
  of~$\mu^+_{\widetilde{T_x}}$ and~$\mu^-_{\widetilde{T_x}}$ . Note
  that~$D_{x,i}$ can be calculated recursively as follows. The main observation
  is that, given a disagreement at~$x$, the random variable
  $|\sigma-\sigma'|_{x,i}$ is the sum of the~$b$ {\em independent} random
  variables $|\sigma-\sigma'|_{z,i-1}$ where~$z$ ranges over the children
  of~$x$. In turn, the random variable $e^{t|\sigma-\sigma'|_{z,i-1}}$ takes the
  value~$D_{z,i-1}$ with probability at most~$\kappa$ (the probability of a
  disagreement at~$z$ given a disagreement at~$x$) and the value~$1$ with the
  remaining probability (since $|\sigma-\sigma'|_{z,i-1}=0$ if there is no
  disagreement at~$z$). Thus, if we let $\delta_i = \max_x D_{x,i}-1$, then
  $\delta_{i+1} \le [1+\kappa\delta_{i}]^b-1 \le e^{\kappa b\delta_{i}}-1 \le
  e^{\alpha\delta_{i}}-1$. We wish to show that, for~$t$ in the above range,
  $\delta_\ell \le 2et\alpha^\ell$, which implies $\Expect_\nu
  \left[e^{t|\sigma-\sigma'|_{r,\ell}}\right] \le \delta_\ell+1 \le
  e^{2et\alpha^\ell}$, as required. In fact, we show by induction that $\delta_i \le
  2t[{\ell+1\over \ell}\cdot\alpha]^i$ for every $\;0\le i\le \ell$. For the
  base case $i=0$, notice that $|\sigma-\sigma'|_{x,0}=1$ when starting from a
  fixed disagreement at~$x$, so $\delta_0=e^t-1 \le 2t$ for $t$ in the given
  range.  For $i+1>0$, we use the fact that 
  $\delta_{i+1} \le e^{\alpha\delta_{i}}-1
  \le {\alpha\delta_{i}\over 1-\alpha\delta_{i}} \le {\ell+1\over \ell}\cdot
  \alpha \delta_{i}$, since by the induction hypothesis 
  $\delta_{i} \le {1\over \alpha (\ell+1)}$ for all $\; 0\le i\le \ell-1$ 
  and $t$ in the given range.
\end{proofof}

\subsection{A crude bound on log-Sobolev via the spectral gap}
\label{subsec:crude}
In this section we state and prove a general bound on~$\csob$ using a bound on
$\cgap$. Although we do not require this bound for the results in this paper,
we believe that it may find applications in the future.
We state the bound for the Ising model, but it can be easily verified that 
it generalizes to any nearest-neighbor spin system on a tree.
  
\begin{theorem}\label{th:crudelogsob}
For the Ising model on the $b$--ary tree, 
$\csob(\mu)=\cgap(\mu)\times\Omega(1/\log n)$. 
In particular, if $\cgap(\mu)=\Omega(1)$ then 
$\csob(\mu)=\Omega(1/\log n)$.
\end{theorem}
It is useful to compare this bound with the well-known bound
$\csob(\mu)=\cgap(\mu)\times \Omega(1/n)$ (see, e.g.,\cite{Saloff}), 
which though much weaker is also more general
(for example, it applies to spin systems on any graph).

Theorem~\ref{th:crudelogsob} is a consequence of the following lemma.
\begin{lemma}
\label{log}
For any $\b$ and~$h$, there exists a constant $c=c(b,\b,h)$ such that, 
for any $x\in T$ and all~$\ell$,
\begin{equation}
  \label{eq:recur}
 \csob(\mu_{T_x}^\tau)^{-1}\le 
\max_{y\prec x,\eta\in\Omega_T^\tau}\{\csob(\mu_{T_y}^\eta)^{-1}\} +
c\cdot\,\cgap(\mu_{T_x}^\tau)^{-1}\,.  
\end{equation}
\end{lemma}
This lemma immediately implies Theorem~\ref{th:crudelogsob}, once
we notice that $\cgap(\mu_{T_x}^\eta)\ge c'\cdot\cgap(\mu_T^\tau)$
for a constant~$c'=c'(b,\beta,h)$ and every $x\in T$ and
$\eta\in \Omega_T^\tau$, as can easily be checked.

\begin{proofofnobox}{Lemma \ref{log}} 
For simplicity and w.l.o.g.\ we will prove the 
recursive inequality~(\ref{eq:recur}) only for $T_x=T$ (the whole tree), 
with root~$r$.  Let $f$ be a non--negative function.
We then write (using the entropy version of (\ref{e:var_dec}))
\begin{equation}
  \label{eq:prel.1}
  \Ent(f) = \mu\bigl[\Ent_{\widetilde T}(f)\bigr] +
  \Ent\bigl[\mu_{\widetilde T}(f)\bigr]\,.
\end{equation}
Using the definition of~$\csob$ we have
\begin{align}
  \mu\bigl[\Ent_{\widetilde T}(f)\bigr] &\le \max_{y\prec
    r,\eta\in\Omega_T^\tau}\{\csob(\mu_{T_y}^\eta)^{-1}\}
  \sum_{x\in \widetilde T}\mu\bigl[\Var_{\{x\}}(\sqrt{f})\bigr] \nonumber \\
  &\le \max_{y\prec r,\eta\in\Omega_T^\tau}\{\csob(\mu_{T_y}^\eta)^{-1}\} \,
  \df(\sqrt{f}\bigr)\,.
 \label{eq:prel.0}
\end{align}
The second term on the r.h.s.\ of~(\ref{eq:prel.1}), being the entropy of a
Bernoulli random variable, is bounded above by 
\begin{eqnarray} 
\Ent\bigl[\mu_{\widetilde T}(f)\bigr]&\le& 
 \alpha\Var\Bigl({\textstyle{\sqrt{\mu_{\widetilde T}(f)}}}\,\Bigr)
 \label{eq:prel.2} \\
&\le& \alpha\Var(\sqrt{f})\nonumber \\  
&\le& \alpha\,\cgap(\mu)^{-1}\,\df\bigl(\sqrt{f}\bigr), 
 \label{eq:prel.3}
\end{eqnarray}
where $\alpha\equiv\alpha(p)$ is a constant that depends on 
$p=\mu(\sigma_r=+)$;
specifically $\alpha(p)=\frac{\log(p/1-p)}{2p-1}$ for
$p\ne1/2$, and $\alpha(1/2)=1/2$ (see~\cite{Saloff}).

Putting together~(\ref{eq:prel.0}) and~(\ref{eq:prel.3}), 
the expression in~(\ref{eq:prel.1}) is bounded above by
$$
\Bigl[\,\max_{y\prec r,\eta\in\Omega_T^\tau}\{\csob(\mu_{T_y}^\eta)^{-1}\} \, +\, \alpha \, \cgap(\mu)^{-1}\,\Bigr]\df\bigl(\sqrt{f}\bigr),
$$
so that from the definition of $\csob$ we have 
$$
 \csob(\mu)^{-1} \le \max_{y\prec r,\eta\in\Omega_T^\tau}\{\csob(\mu_{T_y}^\eta)^{-1}\}+
\alpha\,c_{\rm gap}(\mu)^{-1}.\qquad\square $$
\end{proofofnobox}

%%%%%%%%%%%%%%%%%% EXTENSIONS %%%%%%%%%%%%%%%%%%
%\vskip-0.3in\hbox{}
\section{Extensions to other models}
\label{sec:extensions}
%\vskip-0.05in
As we have already indicated, our techniques extend beyond the
Ising model to general nearest-neighbor interaction models on
trees, including those with hard constraints.  In this final
section we mention some of these extensions.  For a fuller
treatment of this material, the reader is referred to the
companion paper~\cite{companion}.

A ({\it nearest neighbor\/}) {\it spin system\/} on a finite
graph~$G=(V,E)$ is specified by a finite set~$S$ of spin values,
a symmetric pair potential $U:S\times S\to\Rset\cup\{\infty\}$,
and a singleton potential $W:S\to\Rset$.  A configuration
$\sigma\in S^V$ of the system assigns to each vertex (site)
$v\in V$ a spin value $\sigma_v\in S$.  The Gibbs distribution
is given by $$
 \mu(\sigma) \propto \exp\Bigl[-\bigl(\sum\nolimits_{xy\in E}U(\sigma_x,\sigma_y)
             +\sum\nolimits_{x\in V}W(\sigma_x)\bigr)\Bigr].  $$
%where $\beta\ge 0$ plays the role of a generalized inverse temperature.
Thus the Ising model corresponds to the case $S=\{\pm 1\}$, 
and $U(s_1,s_2)=-\beta s_1s_2$, $W(s)=-\beta hs$, where~$\beta$ 
is the inverse temperature and $h$~is the external field.
Note that setting $U(s_1,s_2)=\infty$ corresponds to a 
{\it hard constraint}, i.e., spin values $s_1,s_2$ are 
forbidden to be adjacent.  We denote by~$\Omega$ the set of all
{\it valid\/} spin configurations, i.e., those for which
$\mu(\sigma)>0$.

As for the Ising model, we allow boundary conditions which
fix the spin values of certain sites.  We carry over our
notation from the Ising model: thus, e.g., $\mu_A^\tau$ denotes
the Gibbs distribution on a subset $A\subseteq V$ with boundary
condition~$\tau$ on~$\partial A$.

The (heat-bath) Glauber dynamics extends in the obvious way to general
spin systems.  
We first note that, as the reader may easily check, neither the spatial
mixing conditions in Section~\ref{sec:spat} nor their proofs made any
reference to the details of the Ising model.  All of this material
therefore carries over without modification to general spin systems
on trees.

\begin{theorem}\label{th:ext-spat}
The statements of theorems~\ref{t:spat_fast} and~\ref{t:spat_fast_e} 
hold for general nearest-neighbor spin systems on trees.
\end{theorem}

Likewise, the machinery developed in Sections~\ref{sec:gap} 
and~\ref{sec:logsob}
for verifying the conditions VM and EM also extends to general
models, though the details of the calculations are model-specific.
In particular, Theorems~\ref{t:gap_coup} and~\ref{t:em_coup}
relating $\vcond$ and $\econd$ to the coupling quantities~$\kappa$
and~$\gamma$ of Definition~\ref{d:kp_gm} still hold (with very
minor modifications).  Thus all we need to do 
is to carry out the detailed calculations of~$\kappa$ and~$\gamma$
for the model under consideration.
We now state without proof the results of these calculations 
for several models of interest.  For the proofs, together with
further discussion and extensions, the reader is referred to
the companion paper~\cite{companion}.

\subsection{The hard-core model (independent sets)}
\label{subsec:indsets}
In this model $S=\{0,1\}$, and we refer to a site as {\it occupied\/}
if it has spin value~1, and {\it unoccupied\/} otherwise.
The potentials are $$
   U(1,1)=\infty; \quad U(1,0)=U(0,0)=1;\quad W(1)=L; \quad W(0)=0,  $$
where $L\in\Rset$.  The hard constraint here means that no two adjacent 
sites may be occupied, so $\Omega$ can be identified with the set
of all {\it independent sets\/} in~$G$.  Also, the aggregated
potential of a valid configuration is proportional to the number 
of occupied sites.  Hence the Gibbs distribution takes the simple form  $$
   \mu(\sigma) \propto \lambda^{N(\sigma)},  $$
where $N(\sigma)$ is the number of occupied sites and the 
parameter $\lambda=\exp(-L)> 0$, which controls the density 
of occupation, is referred to as the ``activity.''

The hard-core model on a $b$-ary tree undergoes a phase transition
at a critical activity $\lambda=\lambda_0 = {{b^b}\over{(b-1)^{b+1}}}$
(see, e.g., \cite{Spitzer,Kelly}).  
For $\lambda\le\lambda_0$ there is 
a unique Gibbs measure regardless of the boundary condition 
on the leaves, while for $\lambda>\lambda_0$ there are (at least)
two distinct phases, corresponding to the ``odd'' and ``even''
boundary conditions respectively.  The even boundary condition is
obtained by making the leaves of the tree all occupied if the depth
is even, and all unoccupied otherwise.  The odd boundary condition
is the complement of this.  (These boundary conditions are derived from
the two maximum-density configurations on the infinite tree~$\Tree^b$ 
in which alternate levels --- either odd or even --- are completely occupied.)
For $\lambda>\lambda_0$, the probability of occupation of the root
in the infinite-volume Gibbs measure differs for odd and even boundary 
conditions.  
%It is believed that there is also an ``intermediate'' regime 
%$\lambda_0<\lambda<\lambda_1$ analogous to that for the Ising model,
%\marginpar{\tiny Need to clarify here!!!}
%in which the Gibbs measure, while not unique, is insensitive to
%typical boundary conditions (under a suitable definition of ``typical'');
%however, little is known in rigorous terms about this region.
Relatively little is known about the Glauber dynamics for the hard-core
model on trees, beyond the general result of Luby and
Vigoda~\cite{LV99,V01} which ensures a mixing time of $O(\log n)$ (after
translation to our continuous time setting) when
$\lambda<{2\over{b-1}}$.  This result actually holds for any graph~$G$
of maximum degree~$b+1$.

Our results for the Glauber dynamics in the hard-core model mirror those
given earlier for the Ising model.  First, for sufficiently small
activity~$\lambda$ we show that both $\cgap$ and $\csob$ are uniformly
bounded away from zero for {\it arbitrary\/} boundary conditions.
Second, for {\it even\/} (or, symmetrically, odd) boundary conditions,
we get the same result for {\it all\/} activities~$\lambda$.

\begin{theorem}\label{th:indsets}
For the hard-core model on the $n$-vertex $b$-ary tree with boundary
condition~$\tau$, $\cgap(\mu)$ and $\csob(\mu)$ are $\O(1)$ in both of the
following situations :
\begin{itemize}
\item[{\rm (i)}]  $\tau$ is arbitrary, and $\lambda\le\max\left\{{1\over{\sqrt{b}-1}},\lambda_0\right\}${\rm ;}
\item[{\rm (ii)}] $\tau$ is even (or odd), and $\lambda\ge 0$ is arbitrary.
\end{itemize}
\end{theorem}

Part~(ii) of this theorem is analogous to our earlier result 
for the Ising model with \Plus-boundary and zero external field 
at all temperatures.  This is in line with the intuition that the
even boundary eliminates the only bottleneck in the dynamics.
Part~(i) identifies a region in which the mixing time is insensitive
to the boundary condition.  We would expect this to hold throughout
the low-activity region $\lambda\le\lambda_0$, and indeed, by analogy
with the Ising model, also in some intermediate region beyond this.
Our bound in part~(i) confirms this behavior: note that the quantity
${1\over{\sqrt{b}-1}}$ exceeds~$\lambda_0$ for all $b\ge 5$, and indeed 
for large~$b$ it grows as~${1\over{\sqrt{b}}}$ compared to 
the~${1\over b}$ growth of~$\lambda_0$.  Thus for $b\ge 5$ we 
establish rapid mixing in a region above the critical value~$\lambda_0$.  
To the best of  our knowledge this is the first such result.  
(Note that the result of~\cite{LV99,V01}
mentioned earlier establishes rapid mixing for $\lambda<{2\over{b-1}}$,
which is less than~$\lambda_0$ for all~$b$ and so does not even cover
the whole uniqueness region.)  We should also mention
that our coupling analysis of $\cgap$ in this region has consequences
for the infinite volume Gibbs measure itself, implying that
when $\lambda\le{1\over{\sqrt{b}-1}}$ any
$\mu=\lim_{T\to\infty}\mu_T^\tau$ that is the limit of finite
Gibbs distributions for some boundary configuration~$\tau$ is
{\it extremal}, again a new result.  We elaborate
on these points in the companion paper~\cite{companion}.

\subsection{The antiferromagnetic Potts model (colorings)}
\label{subsec:colorings}
In this model $S=\{1,2,\ldots,q\}$, and the potentials are 
$U(s_1,s_2)=\beta\delta_{s_1,s_2}$, $W(s)=0$.  This is the analog
of the Ising model except that the interactions are {\it antiferromagnetic},
i.e., neighbors with unequal spins are favored.
%$$   \hbox{$U(s_1,s_2) = \begin{cases}
%           -1 & \hbox{\rm if $s_1=s_2$,}\\
%           +1 & \hbox{\rm otherwise;}\end{cases}$};\qquad\qquad 
%   \hbox{$W(s)=0$\ \ $\forall s$}.  $$
The most interesting case of this model is when
$\beta=\infty$ (i.e., zero temperature), which introduces
hard constraints.  Thus if we think of
the $q$ spin values as colors, $\Omega$ is the set of {\it proper
colorings\/} of~$G$, i.e., assignments of colors to vertices
so that no two adjacent vertices receive the same color.  The
Gibbs distribution is uniform over proper colorings.  In this
model it is~$q$ that provides the parameterization.  
For background on the model, see~\cite{BW00}.

For colorings on the $b$-ary tree it is well known that, when
$q\le b+1$, there are multiple Gibbs measures; this follows 
immediately from the existence of ``frozen configurations,''
i.e., colorings in which the color of every internal vertex
is forced by the colors of the leaves (see, e.g.,~\cite{BW00}).
Recently Jonasson~\cite{Jon01} proved that,
as soon as $q\ge b+2$, the Gibbs measure
is unique.  Moreover, it is known that there is again an ``intermediate''
region that includes the value $q=b+1$, in which the Gibbs measure,
while not unique, is insensitive to ``typical'' boundary conditions
(chosen from the free measure); see~\cite{BW00}.

The sharpest result known for the Glauber dynamics on colorings
is due to Vigoda~\cite{V00}, who shows that for arbitrary 
boundary conditions the mixing time
is $O(\log n)$ provided $q>{{11}\over{6}}(b+1)$.  Actually
this result holds for {\it any\/} $n$-vertex graph~$G$ of
maximum degree~$b+1$.\footnote{A recent sequence of papers 
\cite{DF,Mol,Hayes} have
reduced the required number of colors further for general graphs,
under the assumption that the maximum degree is $\Omega(\log n)$;
the current state of the art requires $q\ge(1+\epsilon)(b+1)$
for arbitrarily small $\epsilon>0$~\cite{HV03}.  However, these results
do not apply in our setting where the degree~$b+1$ is fixed.}
Our techniques extend this rapid mixing result all the way down
to the critical value $q\ge b+2$ for which uniqueness holds, 
with arbitrary boundary conditions.  
Again, our result is a consequence of the fact that the associated  
log-Sobolev constant is bounded below by a constant independent of~$n$:
\begin{theorem}\label{th:colorings}
For the colorings model on the $n$-vertex $b$-ary tree with $q\ge b+2$
and arbitrary boundary conditions, both $\cgap(\mu)$ and
$\csob(\mu)$ are~$\O(1)$.
\end{theorem}
%\begin{theorem}\label{th:colorings}
%For the colorings model on the $n$-vertex $b$-ary
%tree the following hold:
%\begin{itemize}
%\item[{\rm (i)}] $\csob(\mu)=\O(1)$ for arbitrary boundary conditions when 
%$q\ge b+2${\rm ;}
%\item[{\rm (ii)}] $\cgap(\mu)=\O(1)$ with free boundary when $q=b+1$.
%\end{itemize}
%\end{theorem}

\subsection{The ferromagnetic Potts model}
\label{subsec:Potts}
Here we have $S=\{1,2,\ldots,q\}$ and potentials 
$U(s_1,s_2)=-\beta\delta_{s_1,s_2}$, $W(s)=0$.
% $$ \hbox{$U(s_1,s_2) = \begin{cases}
%           +1 & \hbox{\rm if $s_1=s_2$,}\\
%           -1 & \hbox{\rm otherwise;}\end{cases}$};\qquad\qquad 
%   \hbox{$W(s)=0$\ \ $\forall s$}.  $$
This is a straightforward generalization of the (ferromagnetic)
Ising model studied earlier in the paper, in which the spin at
each site can take one of $q$ possible values, and the aggregated
potential of any configuration depends on the number of adjacent
pairs of equal spins.  There are no hard constraints.

Qualitatively the behavior of this model is similar to that of
the Ising model, though less is known in precise quantitative terms.
Again there is a phase transition at a critical $\beta=\beta_0$,
which depends on~$b$ and~$q$, so that for $\beta>\beta_0$ (and indeed
for $\beta\ge\beta_0$ when $q>2$) there are
multiple phases.  This value~$\beta_0$ does not in
general have a closed form, but it is known~\cite{Haggstrom} that
$\beta_0<\half\ln({{b+q-1}\over{b-1}})$ for all $q>2$.  (For $q=2$,
this value is exactly~$\beta_0$ for the Ising model as quoted earlier.)

Using our techniques, we are able to prove the following:
\begin{theorem}\label{t:Potts}
For the Potts model on the $n$-vertex $b$-ary tree, $\cgap(\mu)$
and $\csob(\mu)$ are $\Omega(1)$ in all of the following situations:
\begin{itemize}
\vskip-0.2in\hbox{}
\item[{\rm (i)}] the boundary condition is arbitrary and $\b <
  \max\set{\b_0, {1\over 2}\ln({{\sqrt{b}+1}\over {\sqrt{b}-1}})}${\rm ;}

\vskip-0.2in\hbox{}
\item[{\rm (ii)}] the boundary condition is constant (e.g., all sites
  on the boundary have spin~$1$) and~$\b$ is arbitrary{\rm ;}

\vskip-0.2in\hbox{}
\item[{\rm (iii)}] the boundary is free (i.e., the boundary
spins are unconstrained) and $\b < \b_1$,
  where~$\b_1$ is the solution to the equation ${e^{2\b_1}-1\over e^{2\b_1} +q-1}
  \cdot{e^{2\b_1}-1\over e^{2\b_1} +1} = {1\over b}$.
\end{itemize}
\end{theorem}

Part~(i) of this theorem shows that $\cgap$ and $\csob$ are $\Omega(1)$
for arbitrary boundaries throughout the uniqueness region; also, since 
${1\over 2}\ln({{\sqrt{b}+1}\over {\sqrt{b}-1}}) \ge
\half\ln({{b+q-1}\over{b-1}})>\beta_0$
when $q\le 2(\sqrt{b}+1)$, this result extends into the multiple
phase region for many combinations of $b$ and~$q$.
Part~(ii) of the theorem is an analog of our earlier results
for the Ising model with \Plus-boundaries at
all temperatures.  Part~(iii) is of interest for two reasons.
First, since $\beta_1>\beta_0$ always, it exhibits a natural boundary 
condition under which $\cgap$ and $\csob$ are $\Omega(1)$
beyond the uniqueness region 
(but not for arbitrary~$\beta$) {\it for all\/} combinations of~$b$ and~$q$.  
Second, because of an intimate connection between the free boundary
case and so-called ``reconstruction problems'' on trees~\cite{Elchsurvey} 
(in which the edges are noisy channels and the goal is to reconstruct a
value transmitted from the root), we obtain an alternative proof of the
best known value of the noise parameter under which reconstruction is
impossible~\cite{MP01}.
Indeed, a slight strengthening of part~(iii)
allows us to marginally improve on this threshold.
Again, we spell out the details in~\cite{companion}.

%Qualitatively the behavior of this model is similar to that of
%the Ising model, though much less is known in precise quantitative terms.
%In particular, the ``intermediate''
%temperature regime (where the Gibbs measure is not unique but
%typical boundary conditions have no effect) is only very partially
%understood, and the location of the second critical point~$\beta_1$
%is not known.  
%
%Using our techniques, we are able to prove that $\csob(\mu)=\Omega(1)$
%for this model with 
%arbitrary boundary conditions for all $\beta<f(b,q)$, 
%where~$f$ is a certain function of the degree
%of the tree and the number of spin values.  The function~$f$
%is not easy to describe explicitly, so we defer the details
%to~\cite{companion}; however, these seem to be the first
%such results for the Glauber dynamics in this model.  
%Our analysis also has other implications,
%such as a bound on the maximum noise that can be tolerated on the 
%edges of the tree, viewed as noisy communication channels, 
%such that the value transmitted from the root can be reliably
%reconstructed; for a survey on this problem, see~\cite{Elchsurvey}.  
%The threshold that we obtain for this problem is marginally
%better than the best previously known, which is derived in~\cite{MP01}.
%Again, we spell out this application in~\cite{companion}.

%%%%%%%%%%%%%%%%%% APPENDIX %%%%%%%%%%%%%%%%%%
%\newpage
\section{Proofs omitted from the main text}
\label{sec:supplement}
In this final section, we supply the proofs of some technical
lemmas that were omitted from the main text.
%\subsection{Proofs for Section~\protect\ref{sec:spat}}
%\begin{proofof}{Lemma~\ref{l:decay_vcond}}
\subsection{Proof of Lemma~\protect\ref{l:decay_vcond}}
The lemma in fact holds in a more general setting, where in place
of~$\widetilde{T_x}$ and~$B_{x,\ell}$ we think of two arbitrary subsets~$A,B$ 
such that~$A\cup B = T_x$.  Also, in this proof 
we write $\nu=\mu^{\eta}_{T_x}$ and~$\Var$ and~$\Ent$ for variance and entropy
with respect to~$\nu$. For part~(i) we will show that if
for any function~$g$ that does not depend on~$B$ we have
$\Var[\nu_A(g)] \le \varepsilon\cdot \Var(g)$,
then for any function~$f$, $$
\Var[\nu_A(f)] \le {2(1-\varepsilon) \over 1-2\varepsilon}\cdot\nu[\Var_B(f)]+
{2\varepsilon\over 1-2\varepsilon}\cdot \nu[\Var_A(f)].
$$

Notice that by the convexity of variance we have~$\Var(g_1+g_2) \le
2[\Var(g_1)+\Var(g_2)]$ for any two functions~$g_1,g_2$. We therefore write
\begin{eqnarray*}
\Var[\nu_A(f)] &=& \Var[\nu_A(f)-\nu_A(\nu_B(f))+\nu_A(\nu_B(f)]\\
&\le&2\Var[\nu_A(f-\nu_B(f))] + 2\Var[\nu_A(\nu_B(f))]\\
&\le& 2\Var[f-\nu_B(f)] + 2\varepsilon \Var[\nu_B(f)]\\
&=& 2\nu[\Var_B(f)] + 2\varepsilon (\Var[\nu_A(f)] + \nu[\Var_A(f)]-\nu[\Var_B(f)]),
\end{eqnarray*}
where we used the facts that~$\Var[f-\nu_B(f)]=\nu[\Var_B(f)]$ and that
$\Var[\nu_A(f)] + \nu[\Var_A(f)]=\Var[\nu_B(f)] + \nu[\Var_B(f)] = \Var(f)$ as
in~(\ref{e:var_dec}).  We therefore conclude that $\Var[\nu_A(f)] \le
{2(1-\varepsilon)\over 1-2\varepsilon}\cdot \nu[\Var_B(f)] + {2\varepsilon\over
  1-2\varepsilon} \cdot \nu[\Var_A(f)]$, as required.

We proceed to part~(ii). Here we have to show that if for any 
non-negative function~$g$ that does not depend~$B$ we have
$\Ent[\nu_A(g)] \le \varepsilon\cdot \Ent(g)$,
then for any non-negative function~$f$,
\begin{equation}\label{e:econd}
\Ent[\nu_A(f)] \le {1 \over 1-\varepsilon'}\cdot\nu[\Ent_B(f)]+
{\varepsilon'\over 1-\varepsilon'}\cdot \nu[\Ent_A(f)],
\end{equation}
where~$\varepsilon'=\sqrt{\varepsilon}/ p$ and~$p$ stands for the minimum
non-zero probability of any configuration in~$T_x\setminus A$. We will in fact
show that
\begin{equation}\label{e:ent_cesi}
\Ent(f) \le
{1\over 1-\varepsilon'} (\nu[\Ent_A(f)]+\nu[\Ent_B(f)]),
\end{equation}
which implies~(\ref{e:econd}) since~$\Ent[\nu_A(f)]=\Ent(f)-\nu[\Ent_A(f)]$.

Before we go on with the proof, let us review some properties of entropy. First,
by definition, $\Ent(f)=\nu(f\log {f\over \nu (f)})$ and $\nu[\Ent_A(f)]=\nu(f\log
{f\over \nu_A(f)})$. Also, by the variational characterization of entropy 
we have $\nu_A(f\log{g\over \nu_A (g)})\le\Ent_A(f)$ for all non-negative
functions~$f$ and~$g$.

We can now proceed with the proof of~(\ref{e:ent_cesi}) by writing
\begin{eqnarray*}
\Ent(f)&=& \nu\left[f\log {f\over \nu_B(f)}\right] + \nu\left[f\log {\nu_B(f)\over
  \nu_A(\nu_B(f))}\right] + \nu\left[{{f\log {\nu_A(\nu_B(f))}}\over{\nu(f)}}\right]\\
& \le & \nu\left[f\log {f\over \nu_B(f)}\right] + \nu\left[f\log {f\over
  \nu_A(f)}\right] + \nu \left[f\log {\nu_A(\nu_B(f))\over \nu(f)}\right]\\
& = & \nu[\Ent_B(f)] + \nu[\Ent_A(f)] + \nu \left[\nu_A(f)\log {\nu_A(\nu_B(f))\over \nu(f)}\right]\enspace.
\end{eqnarray*}
Therefore,~(\ref{e:ent_cesi}) will follow once we show that~$\nu \left[\nu_A(f)\log
{\nu_A(\nu_B(f))\over\nu(f)}\right] \le \varepsilon'\Ent(f)$. We use the following
claim in order to get this bound.
\begin{claim}\label{c:ent_two}
Let~$\mu$ be a probability measure over a space~$\Omega$ where the probability
 of any~$\sigma\in\Omega$ is either zero or at least~$p$. Then for any two
 non-negative functions~$f$ and~$g$ over~$\Omega$ we have
 $$\mu\left[f\log {g\over \mu g}\right] \le {1\over p}\sqrt{{\mu(f)\over \mu(g)}\cdot
   \Ent(f)\cdot \Ent(g)},$$
 where~$\Ent$ is taken w.r.t. to~$\mu$.
\end{claim}

Assuming Claim~\ref{c:ent_two}, we conclude that
$$\nu\left[\nu_A(f)\log
{\nu_A(\nu_B(f))\over\nu(f)}\right] \le {1\over
  p}\sqrt{\Ent[\nu_A(f)]\cdot\Ent[\nu_A(\nu_B(f))]}
\;\le$$
$${1\over
  p}\sqrt{\varepsilon\cdot\Ent[\nu_A(f)]\cdot\Ent[\nu_B(f)]} \;\le\;
{1\over p}\sqrt{\varepsilon}\Ent(f),$$
completing the proof of Lemma~\ref{l:decay_vcond}. We note that, since
neither~$\nu_A(f)$ nor~$\nu_A(\nu_B(f))$ depends on~$A$, the effective
probability space in the above derivation is the marginal over~$T_x\setminus A$,
so indeed~$p$ can be taken as the minimum marginal
probability of configurations restricted to~$T_x\setminus A$.

It remains to prove claim~\ref{c:ent_two}. Consider two arbitrary
non-negative functions~$f$ and~$g$. 
Let~$\chi$ be the indicator function of the event
that~$g\ge \mu(g)$.  Clearly, $\chi\log {g\over \mu(g)} \ge 0$
while $(1-\chi)\log {g\over \mu(g)} \le 0$. Also, since $\mu\left[\log {g\over
  \mu(g)}\right]\le \log \mu\left[{g\over \mu(g)}\right]=0$ then
$\mu\left[(1-\chi)\log {g\over \mu(g)}\right] \le -\mu\left[\chi\log
  {g\over \mu(g)}\right]$.  Letting~$f_{\mathrm{max}}$ 
and~$f_{\mathrm{min}}$ be the maximum and minimum values of~$f$
respectively over configurations with non-zero probability, we get:
\begin{eqnarray*}
\mu \left[f\log{g\over\mu(g)}\right] & = & \mu\left[\chi f\log
{g\over\mu(g)}\right] + \mu \left[(1-\chi)f\log
{g\over\nu(g)}\right]\\
& \le & f_{\mathrm{max}} \cdot\mu \left[\chi\log{g\over\mu(g)}\right] +
f_{\mathrm{min}} \cdot\mu \left[(1-\chi)\log{g\over\mu(g)}\right]\\
& \le & (f_{\mathrm{max}} - f_{\mathrm{min}})\cdot\mu \left[\chi\log{g\over\mu(g)}\right]\\
& \le & {1\over p}\cdot \|f-\mu(f)\|_1 \cdot\mu\left[\chi\left({g\over\mu(g)}-1\right)\right]\\
& = & {1\over 2p\cdot\mu(g)}\cdot\|f-\mu(f)\|_1 \cdot\|g- \mu(g)\|_1\\
& \le & {1\over p}\sqrt{{\mu(f)\over\mu(g)}\cdot \Ent(f)\cdot\Ent(g)},
\end{eqnarray*}
where we wrote~$\|\cdot\|_1$ for the~$\ell_1$ norm with respect to~$\mu$ and used
 the fact that~$\|f-\mu(f)\|_1^2 \le 2\mu(f)\Ent(f)$ for any non-negative
function~$f$ (see,
e.g.,~\cite{Saloff}). The proof of Claim~\ref{c:ent_two} is now complete.
\stopproof
%\end{proofof}

\subsection{Proof of reverse direction of Theorem~\protect\ref{t:spat_fast_e}} 
%\begin{proofof}{reverse direction of Theorem~\ref{t:spat_fast_e}} 
In the main text we proved the forward direction of
%\marginpar{\tiny READ THIS WHOLE PROOF CAREFULLY!}
Theorem~\ref{t:spat_fast_e}. Here we prove the reverse direction, i.e.,
that $\min_{x,\h}\csob(\mu_{\widetilde T_x}^\h)=\O(1)$
implies~$\econd(\ell,ce^{-\th\ell})$ for all~$\ell$,
where~$c=c(b,\b,h)$ and~$\th=\th(b,\b,h)$ are constants independent
of~$\ell$.  To do this, we follow the same line of reasoning as
in the proof of Theorem~\ref{main_fabio}: namely, we establish
the strong concentration property of the functions~$g_s^{(\ell)}$
as in Section~\ref{subsec:econd_conc} and then appeal to 
Theorem~\ref{th:conc_econd}.
The proof of concentration is accomplished via hypercontractivity
bounds, assuming the above condition on~$\csob$.

For a function~$f$, let~$\Lambda_f\subseteq T$ denote the subset of
sites on whose spins~$f$ depends. We then have:
\begin{lemma}\label{l:sob_norm}
Let $\nu$ be any Gibbs measure on~$T$, $f$~any function, and $B$
any subset that includes all sites within distance~$\ell$ from~$\Lambda_f$.
Then there exists a constant~$\th'$, depending only on the degree~$b$,
such that $$
  \|\nu_B(f)-\nu(f)\|_q \;\le\; 3e^{-\csob(\nu)\th'\ell}|\Lambda_f|\|f-\nu(f)\|_{\infty}\,,$$
where $q=1+e^{\csob(\nu)\th'\ell}$ and norms are taken w.r.t.~$\nu$.
\end{lemma}
We first assume Lemma~\ref{l:sob_norm} and complete the proof of the
reverse direction of Theorem~\ref{t:spat_fast_e}.

For simplicity, we verify~$\econd$ only for the case $T_x=T$ (the whole
tree), with root~$r$.
Recall the functions~$g_s^{(\ell)}$ from Section~\ref{subsec:econd_conc},
the fact that~$g_s^{(\ell)}=\mu_{B_{r,\ell}}(g_s)$ by definition, and that~$g_s$
depends only on the spin at~$r$.  Applying Lemma~\ref{l:sob_norm} with
$\nu=\mu$, $f=g_s$, and $B=B_{r,\ell}$, together with the fact
that~$\csob(\mu)=\Omega(1)$ by hypothesis, we conclude that there exists a
constant~$\th''$ such that
$$
\|g_s^{(\ell)}-1\|_q \le
3e^{-\th''\ell}\ninf{g_s - 1} \le 3e^{-\th''\ell}/p_{\min},
$$ 
where~$q=1+e^{\th''\ell}$ and norms are taken
w.r.t.~$\mu$.  Therefore, using a Markov inequality, there
exist constants~$\ell_0$ and $\th$ such that, for all~$\ell\ge \ell_0$,
$$
\mu^+_{\widetilde{T}}\bigl[|g_s^{(\ell)} - 1| > e^{-\th\ell}\bigr] \;\le\;
e^{-2e^{\th \ell}}.
$$
This establishes the strong concentration property of~$g_s^{(\ell)}$
as in~(\ref{eq:strconc}), from which $\econd$ follows by 
Theorem~\ref{th:conc_econd}.
\stopproof
%\end{proofof}
\begin{remark*}
A similar claim to Lemma~\ref{l:sob_norm} was proved in~\cite{SZ2} in
the context of~$\Z^d$; we reprove it below for completeness. The
proof, as well as the fact that a $\O(1)$ logarithmic Sobolev constant
implies $\econd(\ell,ce^{-\th\ell})$, applies to general, finite range
models on any graph of bounded degree.
\end{remark*}
\begin{proofof}{Lemma~\ref{l:sob_norm}}
The proof has two main ingredients: the first
is a bound on the speed at which information propagates under the
Glauber dynamics, while the second is a standard relationship
between~$\csob$ and hypercontractivity bounds.

Let~$P_t=\nep{t\cL}$ stand for the transition kernel at time~$t$ (as
discussed in Section~\ref{sec:prelims}) of the dynamics under
consideration, reversible w.r.t.\ the Gibbs measure~$\nu$,
and let~$P_t^{B}$ stand for the transition kernel of a
modified dynamics where the spins of the sites outside the subset~$B$
are fixed to their values at time zero (the sites inside~$B$ being
updated according to the same rule as in the original dynamics). It is
well known (see, e.g.,~\cite{SZ2}) that there exists a constant~$k_0$
depending only on~$b$ (or on the degree of the graph in the general case)
and the maximum flip rate
$\max_x \ninf{c_x}$ (which is bounded by~1 in the case of the
heat bath dynamics) such that, for any function~$f$, any~$t$ and any
subset~$B$ that includes all sites within distance~$k_0 t$
of~$\Lambda_f$,
\begin{equation}\label{e:spd_info}
\ninf{P_t f-P_t^B f} \le 2e^{-t}|\Lambda_f|\ninf{f}.
\end{equation}
Equation~(\ref{e:spd_info}) is a manifestation of the fact that it takes
at least~${\ell\over k_0}$ time before the spin at a site can become 
sensitive to the configuration at distance~$\ell$ from it.

The second ingredient we need is a hypercontractivity bound. From Gross's
integration lemma (see, e.g.,~\cite{bible}), we have 
$\|P_tf\|_q \le \|f\|_2$ for any mean-zero function~$f$, any~$t$, 
and~$2\le q\le 1+e^{\csob t}$, where $\csob=\csob(\nu)$.  Adding
to this the fact that~$\cgap>\csob$, we may write
\begin{equation}\label{e:hypc}
   \|P_t f\|_q = \|P_{t/2}(P_{t/2} f)\|_q \le \| P_{t/2} f\|_2
         \le e^{-\cgap t/2}\|f\|_2 \le e^{-\csob t/2}\|f\|_2,
\end{equation}
where~$q=1+e^{\csob t/2}$ and we used the fact that~$\cgap$ bounds the
rate of decay of the~$L^2$ norm.

We now conclude the proof of Lemma~\ref{l:sob_norm} as follows. Without loss of
generality, consider an arbitrary function~$f$ with~$\nu(f)=0$. Let~$\ell$ be
arbitrary, and~$B$ be a subset that includes all sites within distance~$\ell$
of~$\Lambda_f$. Then, for $t=\ell/k_0$ and~$q=1+e^{\csob t/2}$, we have
\begin{eqnarray*}
\|\nu_B(f)\|_q & = & \|\nu_B(P_t^Bf)\|_q\\
& \le & \|P_t^B f\|_q\\
& \le & \|P_t^Bf - P_tf\|_q + \|P_tf\|_q\\
& \le & 2e^{-t}|\Lambda_f|\ninf{f} + e^{-\csob t/2}\|f\|_2\\
& \le & 3e^{-\csob \th \ell}|\Lambda_f|\ninf{f}\,,
\end{eqnarray*}
taking the constant $\th=1/2k_0$ (and using the fact that $\csob\le 1$).
\end{proofof}

%\subsection{Proofs for Section~\protect\ref{sec:logsob}}
%\begin{proofof}{Lemma~\ref{En_in_gen}}
\subsection{Proof of Lemma~\protect\ref{En_in_gen}}
We split our analysis of $\nu(f_1f_2)^2$ into three cases:
\begin{enumerate}[(a)]
\item $\Ent_\nu(f_2) \ge \ov\d$;

\smallno
\item $\d< \Ent_\nu(f_2)< \ov \d$;

\smallno
\item $\Ent_\nu(f_2)\le \d$.
\end{enumerate}

\smallno
\emph{Case (a)}. We simply bound
\begin{equation*}
 \nu(f_1f_2)^2\le \ninf{f_1}^2\nu(f_2)^2\le
1\le 
 \d\,\Ent_{\nu}(f_2)\,.
\end{equation*}
\emph{Case (b)}. We use the  
\emph{entropy inequality} (see, e.g.,~\cite{bible}),
which states that for any $t>0$,
\begin{equation}
  \label{entro_in}
  \nu(f_1f_2) \le \frac{1}{t}\log \nu(\nep{tf_1}) + \frac{1}{t}\Ent_\nu(f_2)\,.
\end{equation}
We
choose the free parameter $t$ in (\ref{entro_in}) equal to  
$\sqrt{\Ent_\nu(f_2)/\d}$. Notice that, by construction, $1<t <\d^{-1}$.
Using the assumption $\nu(|f_1|> \d)\le \nep{-2/\d}$ together with
$\ninf{f_1}\le 1$, we get
\begin{align*}
\nu(f_1f_2)^2 &\le \Bigl[ \ov t \log\bigl(\nep{t\d} + \nep{t-
 2/\d}\bigr)+ \sqrt{\d \Ent_\nu(f_2)}\ \Bigr]^2 \\
&\le \Bigl[ c_1\, \d + \sqrt{\d \Ent_\nu(f_2)}\ \Bigr]^2 \le c_2\,\d \Ent_\nu(f_2)
\end{align*}
for suitable  numerical constants $c_1,c_2$.

\medno
\emph{Case (c)}. Again we use the entropy inequality with 
$t=\sqrt{\Ent_\nu(f_2)/\d} \le 1$, but we now simply bound the Laplace
transform $\nu(\nep{tf_1})$ by a Taylor
expansion (in $t$) up to second order:
\begin{align*}
\ov t \log\nu(\nep{tf_1}) &\le \ov t \log\Bigl( 
1 + \nep{}\frac{t^2}{2} \nu(f_1^2) \Bigr)
\le \nep{}\frac{t}{2} \bigl[\,\d^2 + \nep{-2/\d}\,\bigr] \\
&= \ov 2\nep{}\bigl[\,\d^2 + \nep{-2/\d}\,\bigr]\sqrt{\Ent_\nu(f_2)/\d},
\label{th.11}
\end{align*}
which by~(\ref{entro_in}) implies 
\begin{equation*}
  \nu(f_1f_2)^2 \le \Bigl[\,\frac{e}{2 \sqrt{\d}}\,\bigl(\,\d^2 + \nep{-2/\d}\,\bigr)
  + \sqrt{\d}\, \Bigr]^2 \Ent_\nu(f_2) \le c_3 \,\d \Ent_\nu(f_2) 
\end{equation*}
for another numerical constant $c_3$.
\stopproof
%\end{proofof}

%%%%%%%%%%%%%%%%%% REFERENCES %%%%%%%%%%%%%%%%%%
\section*{References}
\begin{list}%
        {[\arabic{refcount}]}{\usecounter{refcount}
         \setlength{\itemsep}{0pt}\setlength{\leftmargin}{0.8truecm}}

% \bibitem{Aldous}
% {\sc D.~Aldous},
% ``Random walks on finite groups and rapidly mixing Markov chains,''
% {\it S\'eminaire de Probabilites XVII},
% Lecture Notes in Mathematics~{\bf 986}, pp.~243--297,
% Springer, Berlin, 1983. 

\bibitem{bible} 
{\sc C.~An{\'e}, S.~Blach{\`e}re, D.~Chafa{\"\i}, P.~Foug{\`e}res, 
I.~Gentil, F.~Malrieu, C.~Roberto} and {\sc G.~Scheffer},
``Sur les in\'egalit\'es de {S}obolev logarithmiques,'' 
{\it Soci\'et\'e Math\'ematique de France}, 2000.

\bibitem{Baxter}
{\sc R.J.~Baxter},
{\it Exactly solved models in statistical mechanics},
Academic Press, London, 1982.

\bibitem{BKMP} 
{\sc N.~Berger, C.~Kenyon, E.~Mossel} and {\sc Y.~Peres},
``Glauber dynamics on trees and hyperbolic graphs,'' preprint (2003).
Preliminary version:
{\sc C.~Kenyon, E.~Mossel} and {\sc Y.~Peres},
``Glauber dynamics on trees and hyperbolic graphs,''
{\it Proc.\ 42nd IEEE Symposium on Foundations of Computer Science\/}~(2001),
pp.~568--578.

\bibitem{Cesi'}
{\sc L.~Bertini, N.~Cancrini} and {\sc F.~Cesi},
``The spectral gap for a Glauber-type dynamics in a continuous gas,''
{\it Ann.\ Inst.\ H.~Poincar{\'e} Probab. Statist.}~{\bf 38} (2002),
pp.~91--108.

\bibitem{BRSSZ}
{\sc P.~Bleher, J.~Ruiz, R.H.~Schonmann, S.~Shlosman} and {\sc V.~Zagrebnov},
``Rigidity of the critical phases on a Cayley tree,''
{\it Moscow Mathematical Journal\/}~{\bf 1} (2001), pp.~345--363.

\bibitem{BRZ}
{\sc P.~Bleher, J.~Ruiz} and {\sc V.~Zagrebnov},
``On the purity of the limiting Gibbs state for the Ising model
on the Bethe lattice,''
{\it Journal of Statistical Physics\/}~{\bf 79} (1995), pp.~473--482.

\bibitem{BM}
{\sc T.~Bodineau} and {\sc F.~Martinelli},
``Some new results on the kinetic Ising model in a pure phase,''
{\it Journal of Statistical Physics\/}~{\bf 109}~(1), 2002.

\bibitem{BW00}
{\sc G.~Brightwell} and {\sc P.~Winkler},
``Random colorings of a Cayley tree,''
{\it Contemporary combinatorics}, Bolyai Society Mathematical Studies~{\bf 10},
J\'anos Bolyai Math. Soc., Budapest, 2002, pp.~247--276.

\bibitem{Cesi} 
{\sc F.~Cesi},
``Quasi-factorization of the entropy and logarithmic Sobolev inequalities for
Gibbs random fields,''
{\it Probability Theory and Related Fields\/}~{\bf 120} (2001), 
pp.~569--584.

\bibitem{CCST} 
{\sc J. T.~Chayes, L.~Chayes, J.P.~Sethna} and {\sc D.J.~Thouless},
``A mean field spin glass with short-range interactions''
{\it Communications in Mathematical Physics\/}~{\bf 106} (1986),
pp.~41--89.

\bibitem{DKS}
{\sc R.~Dobrushin, R.~Koteck\'y} and {\sc S.~Shlosman},
``Wulff Construction. A Global Shape From Local
   Interaction'' {\it Translation of Math. Monographs\/}, AMS {\bf 104} (1992). 
%\bibitem{DS-C}
%{\sc P.~Diaconis} and {\sc L.~Saloff-Coste},
%``Comparison theorems for reversible Markov chains,''
%{\it Annals of Applied Probability\/}~{\bf 3} (1993), pp.~696--730.

\bibitem{DF}
{\sc M.~Dyer} and {\sc A.~Frieze},
``Randomly colouring graphs with lower bounds on girth and maximum degree,''
{\it Proceedings of the 42nd Annual IEEE Symposium on Foundations of
Computer Science}, 2001, pp.~579--587.

\bibitem{EKPS}
{\sc W.~Evans, C.~Kenyon, Y.~Peres} and {\sc L.J.~Schulman},
``Broadcasting on trees and the Ising model,''
{\it Annals of Applied Probability\/}~{\bf 10} (2000), pp.~410--433.

\bibitem{FH}
{\sc D.~Fisher} and {\sc D.~Huse},
``Dynamics of droplet fluctuations in pure and random Ising systems,''
{\it Physics Review~B\/}~{\bf 35}~(13), 1987.

% \bibitem{FK}
% {\sc A.~Frieze} and {\sc R.~Kannan},
% ``Log-Sobolev inequalities and sampling from log-concave distributions,''
% {\it Annals of Applied Probability\/}~{\bf 9} (1999), pp.~14--26.

\bibitem{Georgii}
{\sc H.-O. Georgii},
{\it Gibbs measures and phase transitions}, 
de Gruyter Studies in Mathematics~{\bf 9}, Walter de Gruyter \& Co., 
Berlin, 1988.

\bibitem{Haggstrom}
{\sc O.~H\"aggstr\"om},
``The random-cluster model on a homogeneous tree,''
{\it Probability Theory and Related Fields\/}~{\bf 104} (1996),
pp.~231--253.

\bibitem{Hayes}
{\sc T.P.~Hayes},
``Randomly coloring graphs with girth five,''
{\it Proceedings of the 35th Annual ACM Symposium on Theory of Computing}, 
2003, pp.~269--278.

\bibitem{HV03}
{\sc T.P.~Hayes} and {\sc E.~Vigoda},
``A non-Markovian coupling for randomly sampling colorings,''
to appear in {\it Proceedings of the 44th Annual IEEE Symposium
on Foundations of Computer Science}, 2003.

\bibitem{Ioffe}
{\sc D.~Ioffe},
``A note on the extremality of the disordered state for the Ising model
on the Bethe lattice,''
{\it Letters in Mathematical Physics\/}~{\bf 37} (1996), pp.~137--143.

\bibitem{Ioffe2}
{\sc D.~Ioffe},
``Extremality of the disordered state for the Ising model on general trees,''
{\it Progress in Probability\/}~{\bf 40} (1996), pp.~3--14.

% \bibitem{Ising}
% {\sc E.~Ising},
% ``Beitrag zur Theorie des Ferromagnetismus,''
% {\it Zeitschrift f\"ur Physik\/}~{\bf 31} (1925), pp.~253--258.

% \bibitem{JS}
% {\sc M.~Jerrum} and {\sc J.-B.~Son},
% ``Spectral gap and log-Sobolev constant for balanced matroids,''
% {\it Proceedings of the 43rd IEEE Symposium on Foundations of 
% Computer Science}, 2002, pp.~721--729.

% \bibitem{JSTV}
% {\sc M.~Jerrum, J.-B.~Son, P.~Tetali} and {\sc E.~Vigoda},
% ``Elementary bounds on Poincar\'e and log-Sobolev constants for 
% decomposable Markov chains,''
% Preprint, 2003.

\bibitem{Jon01}
{\sc J.~Jonasson}
``Uniqueness of uniform random colorings of regular trees,''
{\it Statistics \& Probability Letters}~{\bf 57} (2002), pp.~243--248.

\bibitem{JS99}
{\sc J.~Jonasson} and {\sc J.E.~Steif},
``Amenability and phase transition in the Ising model,''
{\it Journal of Theoretical Probability\/}~{\bf 12} (1999), pp.~549--559.

\bibitem{Kelly}
{\sc F.P.~Kelly},
``Stochastic models of computer communication systems,''
{\it Journal of the Royal Statistical Society\/}~B~{\bf 47} (1985),
pp.~379--395.

% \bibitem{KMP} 
% {\sc C.~Kenyon, E.~Mossel} and {\sc Y.~Peres},
% ``Glauber dynamics on trees and hyperbolic graphs,''
% {\it Proc. 42nd IEEE Symposium on Foundations of Computer Science\/}~(2001),
% pp.~568--578.

\bibitem{L}
{\sc R.~Lyons},
{\it Phase transitions on non amenable graphs},
J.Math.Phys~{\bf 41}, pp.~1099--1127, 2000.

% \bibitem{Latala} 
% {\sc R.~Latala}, ``Between Sobolev and Poincar\'e,'' 
% {\it Geometric aspects of functional analysis}, {\it Lecture Notes in
% Mathematics}~{\bf 1745}, pp.~147--168, Springer, Berlin, 2000.

\bibitem{Ledoux} 
{\sc M.~Ledoux}, ``The concentration of measure phenomenon,'' 
{\it Mathematical Surveys and Monographs}~{\bf 89}, American 
Mathematical Society, 1981.

% \bibitem{Lenz}
% {\sc W.~Lenz},
% ``Beitrag zum Verst\"andnis der magnetischen Erscheinungen in festen 
% K\"orpern,''
% {\it Zeitschrift f\"ur Physik\/}~{\bf 21} (1920), pp.~613--615.

\bibitem{Yau} 
{\sc S.L.~Lu} and {\sc H.T.~Yau},
``Spectral gap and logarithmic Sobolev inequality for Kawasaki and
Glauber dynamics''
{\it Comm. Math. Phys.}~{\bf 156} (1993), pp.~399--433.

\bibitem{LV99}
{\sc M.~Luby} and {\sc E.~Vigoda},
``Fast convergence of the Glauber dynamics for sampling independent sets,''
{\it Random Structures \& Algorithms\/}~{\bf 15} (1999), pp.~229--241.

\bibitem{Martinelli} 
{\sc F.~Martinelli},
``Lectures on Glauber dynamics for discrete spin models,''
{\it Lectures on Probability Theory and Statistics (Saint-Flour, 1997)},
Lecture notes in Mathematics~{\bf 1717}, pp.~93--191,
Springer, Berlin, 1998.

\bibitem{MO1}
{\sc F.~Martinelli} and {\sc E.~Olivieri},
``Approach to equilibrium of Glauber dynamics in the one phase region I: The
attractive case,'' 
{\it Comm. Math. Phys.}~{\bf 161} (1994), pp.~447--486.

\bibitem{MOS}
{\sc F.~Martinelli, E.~Olivieri} and {\sc R.~Schonmann},
``For 2-D lattice spin systems weak mixing implies strong mixing,''
{\it Comm.\ Math.\ Phys.}~{\bf 165} (1994), pp.~33--47.

\bibitem{companion}
{\sc F.~Martinelli, A.~Sinclair} and {\sc D.~Weitz},
``Fast mixing for independent sets, colorings and other models on trees,''
preprint, 2003.

\bibitem{Mol}
{\sc M.~Molloy},
``The Glauber dynamics on colorings of a graph with high girth and
maximum degree,''
{\it Proceedings of the 34th Annual ACM Symposium on Theory of Computing},
2002, pp.~91--98.

\bibitem{Elchsurvey}
{\sc E.~Mossel},
``Survey: Information flow on trees,''
Preprint, October 2002, to appear in DIMACS volume
{\it Graphs, Morphisms and Statistical Physics}.

\bibitem{MP01}
{\sc E.~Mossel} and {\sc Y.~Peres},
``Information flow on trees,''
{\it Annals of Applied Probability}, 2003, to appear.

\bibitem{PW}
{\sc Y.~Peres} and {\sc P.~Winkler},
personal communication.

% \bibitem{Preston}
% {\sc C.J.~Preston},
% {\it Gibbs states on countable sets},
% Cambridge Tracts in Mathematics~{\bf 68}, Cambridge University
% Press, London, 1974.

\bibitem{Saloff}
{\sc L.~Saloff-Coste},
``Lectures on finite Markov chains,''
{\it Lectures on probability theory and statistics (Saint-Flour, 1996)}, 
Lecture notes in Mathematics~{\bf 1665}, pp.~301--413, 
Springer, Berlin, 1997. 

\bibitem{ST98}
{\sc R.H.~Schonmann} and {\sc N.I.~Tanaka},
``Lack of monotonicity in ferromagnetic Ising model phase diagrams,''
{\it Annals of Applied Probability\/}~{\bf 8} (1998), pp.~234--245.

\bibitem{Simon}
{\sc B.~Simon},
{\it The statistical mechanics of lattice gases, Vol.~I},
Princeton Series in Physics, Princeton University Press, 
Princeton, NJ, 1993.

\bibitem{Spitzer}
{\sc F.~Spitzer},
``Markov random fields on an infinite tree,''
{\it Annals of Probability\/}~{\bf 3} (1975), pp.~387--398.

\bibitem{SZ} 
{\sc D.W.~Stroock} and {\sc B.~Zegarlinski},
``The logarithmic Sobolev inequality for discrete spin systems on a lattice,''
{\it Comm. Math. Phys.}~{\bf 149} (1992), pp.~175--194.

\bibitem{SZ2} 
{\sc D.W.~Stroock} and {\sc B.~Zegarlinski},
``On the ergodic properties of Glauber dynamics,''
{\it J. Statist. Phys.}~{\bf 81} (1995), pp.~1007--1019.

\bibitem{V00}
{\sc E.~Vigoda},
``Improved bounds for sampling colorings,''
{\it  Journal of Mathematical Physics\/}~{\bf 41} (2000), pp.~1555--1569.

\bibitem{V01}
{\sc E.~Vigoda},
``A note on the Glauber dynamics for sampling independent sets,''
Electronic Journal of Combinatorics, Volume~8(1), 2001.

\end{list}

\end{document}